\documentclass[authoryear,11pt,preprint]{elsarticle}
\usepackage{color}
\usepackage{lineno}
\usepackage{pdflscape}
\usepackage{amsmath}
\usepackage{booktabs}
\usepackage[hyphens]{url}
\usepackage{hyperref}
\usepackage{multirow}
\usepackage{addlines} %A
\usepackage{adjustbox}
\usepackage{siunitx}
\usepackage{url}
\usepackage{natbib} %When changing reference style
\usepackage[english]{babel}
\usepackage[dvipsnames]{xcolor}
\usepackage{soul}
\usepackage{adjustbox}
\usepackage{verbatim}
\usepackage[margin=2.5cm]{geometry}% by courtesy of Mico

\usepackage{mathtools}
\usepackage{tikz}
\usepackage{pgfplots}

\definecolor{darkblue}{rgb}{0,0.2,0.4}
\usepackage{graphicx}
\usepackage{subcaption}
\usepackage{amssymb}
\usepackage{amsthm}
\usepackage{etoolbox}
\usepackage{algorithm}
\usepackage{algpseudocode}
\usepackage{lipsum}

\usepackage{enumitem}
\usepackage{calc}

\makeatletter
\newcommand*{\algrule}[1][\algorithmicindent]{%
  \hspace*{.2em}% <------------- This is where the rule starts from
  \vrule %height .75\baselineskip depth .25\baselineskip
  \hspace*{\dimexpr#1-.2em-.4pt}%
}

\newcommand{\StatePar}[1]{%
  \State\parbox[t]{\dimexpr\linewidth-\ALG@thistlm}{\strut #1\strut}%
}
\renewcommand{\ALG@beginalgorithmic}{\offinterlineskip}% Remove all interline skips

\newcount\ALG@printindent@tempcnta
\def\ALG@printindent{%
  \ifnum \theALG@nested > 0% is there anything to print
    \ifx\ALG@text\ALG@x@notext% is this an end group without any text?
    \else
      \unskip
      \ALG@printindent@tempcnta=1
      \loop
        \algrule[\csname ALG@ind@\the\ALG@printindent@tempcnta\endcsname]%
        \advance \ALG@printindent@tempcnta 1
        \ifnum \ALG@printindent@tempcnta<\numexpr\theALG@nested+1\relax
      \repeat
        \fi
    \fi
}
\patchcmd{\ALG@doentity}{\noindent\hskip\ALG@tlm}{\ALG@printindent}{}{\errmessage{failed to patch}}
\makeatother

\algrenewcommand\algorithmicend{\strut\textbf{end}}
\algrenewcommand\algorithmicdo{\strut\textbf{do}}
\algrenewcommand\algorithmicwhile{\strut\textbf{while}}
\algrenewcommand\algorithmicfor{\strut\textbf{for}}
\algrenewcommand\algorithmicforall{\strut\textbf{for all}}
\algrenewcommand\algorithmicloop{\strut\textbf{loop}}
\algrenewcommand\algorithmicrepeat{\strut\textbf{repeat}}
\algrenewcommand\algorithmicuntil{\strut\textbf{until}}
\algrenewcommand\algorithmicprocedure{\strut\textbf{procedure}}
\algrenewcommand\algorithmicfunction{\strut\textbf{function}}
\algrenewcommand\algorithmicif{\strut\textbf{if}}
\algrenewcommand\algorithmicthen{\strut\textbf{then}}
\algrenewcommand\algorithmicelse{\strut\textbf{else}}
\algrenewcommand\algorithmicrequire{\strut\textbf{Input:}}
\algrenewcommand\algorithmicensure{\strut\textbf{Output:}}
\let\oldState\State
\renewcommand{\State}{\oldState\strut}

\usepackage{nomencl}
\renewcommand\nomgroup[1]{%
  \item[\bfseries
  \ifstrequal{#1}{A}{Investment planning model sets}{%
  \ifstrequal{#1}{B}{Operational model sets}{%
  \ifstrequal{#1}{C}{Investment planning model parameters}{%
  \ifstrequal{#1}{D}{Operational model parameters}{%
  \ifstrequal{#1}{E}{Investment planning model variables}{%
  \ifstrequal{#1}{F}{Operational model variables}{
\ifstrequal{#1}{G}{Function}{
  }}}}}}}%
]}
\makenomenclature

\usepackage{framed} % Framing content
\usepackage{multicol} % Multiple columns environment
\usepackage{eurosym}
\usepackage[utf8]{inputenc}
\usepackage{newunicodechar}

\DeclareUnicodeCharacter{0229}{\c{e}}
\usepackage[T1]{fontenc}
\usepackage{textcomp}
\usepackage{lmodern}

\usepackage{ltxcmds}

\usepackage{siunitx}
\DeclareSIUnit{\euro}{\texteuro}

\newcommand*{\sieuro}[2][]{\SI[{mode=text,#1}]{#2}{\euro}}

\makeatletter
\newcommand*{\EuroMacro}{}
\protected\def\EuroMacro{%
  \ltx@ifnextchar@nospace\bgroup\sieuro{%
    \ltx@ifnextchar[\sieuro\texteuro
  }%
}

\makeatother

\definecolor{royalpurple}{rgb}{0.58, 0.44, 0.86}

\allowdisplaybreaks
\newcommand\norm[1]{\left\lVert#1\right\rVert}

\algnewcommand{\IIf}[1]{\State\algorithmicif\ #1\ \algorithmicthen}
\algnewcommand{\EndIIf}{\unskip\ \algorithmicend\ \algorithmicif}

\usepackage{setspace}
\hypersetup{colorlinks = true,linkcolor = blue,anchorcolor =red,citecolor = blue,filecolor = red,urlcolor = red, pdfauthor=author}

\journal{}

\hypersetup{ colorlinks=true,
linkcolor=darkblue,
filecolor=darkblue,
citecolor=darkblue,
urlcolor=darkblue
}
\pgfplotsset{compat=1.17}
\begin{document}
\begin{frontmatter}


\title{A stabilised Benders decomposition with adaptive oracles applied to investment planning of multi-region power systems with short-term and long-term uncertainty}


\author[a]{Hongyu Zhang\corref{cor1}}
\ead{hongyu.zhang@ntnu.no}
\author[c]{Nicolò Mazzi}
\ead{mazzi.nicolo@gmail.com}
\author[b]{Ken McKinnon}
\ead{K.McKinnon@ed.ac.uk}
\author[b]{Rodrigo Garcia Nava}
\ead{Rodrigo.Garciana@ed.ac.uk}
\author[a]{Asgeir Tomasgard}
\ead{asgeir.tomasgard@ntnu.no}

\cortext[cor1]{Corresponding author}

\address[a]{Department of Industrial Economics and Technology Management, Norwegian University of Science and Technology, Høgskoleringen 1, 7491, Trondheim, Norway}
\address[b]{School of Mathematics, University of Edinburgh, James Clerk Maxwell Building, Peter Guthrie Tait Road, Edinburgh, EH9 3FD, United Kingdom}
\address[c]{aHead-Research, Corso Svizzera 185, Torino, 10149, Italy}

\begin{abstract}
Benders decomposition with adaptive oracles was proposed to solve large-scale optimisation problems with a column bounded block-diagonal structure, where subproblems differ on the right-hand side and cost coefficients. Adaptive Benders reduces computational effort significantly by iteratively building inexact cutting planes and valid upper and lower bounds. However, Adaptive Benders and standard Benders may suffer severe oscillation when solving a multi-region investment planning problem. Therefore, we propose stabilising Adaptive Benders with the level set method and adaptively selecting the subproblems to solve per iteration for more accurate information. Furthermore, we propose a dynamic level set method to improve the robustness of stabilised Adaptive Benders by adjusting the level set per iteration. We compare stabilised Adaptive Benders with the unstabilised versions of Adaptive Benders with one subproblem solved per iteration and standard Benders on a multi-region long-term power system investment planning problem with short-term and long-term uncertainty. The problem is formulated as multi-horizon stochastic programming. Four algorithms were implemented to solve linear programming with up to $1$ billion variables and $4.5$ billion constraints. The computational results show that: a) for a $1.00\%$ convergence tolerance, the proposed stabilised method is up to $113.7$ times faster than standard Benders and $2.14$ times faster than unstabilised Adaptive Benders ; b) for a $0.10\%$ convergence tolerance, the proposed stabilised method is up to $45.5$ times faster than standard Benders  and unstabilised Adaptive Benders cannot solve the largest instance to convergence tolerance due to severe oscillation and c) dynamic level set method makes stabilisation more robust.
\end{abstract}

\begin{keyword}
Large scale optimisation \sep Multi-stage stochastic programming \sep Multi-horizon stochastic programming \sep Stabilised Benders decomposition with adaptive oracles \sep Level set method

% MSC codes here, in the form: \MSC code \sep code
% or \MSC[2008] code \sep code (2000 is the default)

\end{keyword}
%Highlights
% Presemt a new vision on creating vir
\end{frontmatter}

% \linenumbers

\section{Introduction}
\label{sec:introduction}
Power system infrastructure planning is crucial during the energy transition towards zero emission by $2050$. Optimisation models are widely used for the investment and operational planning of systems. To gain enough environmental and economic insights from such models, sometimes a large-scale problem needs to be modelled, such as \citep{Li2022,Zhang2021}. An investment planning problem can involve many technologies and regions and span over decades with multiple investment periods \citep{Conejo2016}. This can lead to a large-scale optimisation problem that is intractable. Furthermore, investment planning of a power system often faces uncertainty from two time horizons \citep{Kaut2014,lara2020}: a) the uncertainty from the operational time horizon, such as the availability of renewable energy. The operational uncertainty becomes even more crucial for a system with higher penetration of intermittent renewable energy, and b) the uncertainty from the strategic time horizon, e.g., CO$_2$ tax and CO$_2$ budget. Stochastic programming is often used to model uncertainty. However, including uncertainty from both time horizons using multi-stage stochastic programming may lead to a large scenario tree and an intractable model. Most studies on power system investment planning in a multi-horizon framework only consider short-term uncertainty, such as \citep{Backe2022}, partly because of the tractability of the problem. However, short-term and long-term uncertainty can play a decisive role in investment planning. Although there are examples including short-term and long-term uncertainty in a multi-horizon model \citep{hellemo2013multi}, the computational difficulty is not sufficiently addressed. Therefore, we aim to address the computational difficulties of long-term planning problems with short-term and long-term uncertainty. One possible way is to reduce the problem size by using a different modelling approach called multi-horizon stochastic programming \citep{Kaut2014}. Although multi-horizon stochastic programming can reduce the scenario tree significantly, it is essentially multi-stage stochastic programming once both short-term and long-term uncertainty is added and can be intractable when the problem gets large. Another way is to develop an algorithm that can efficiently solve a class of large-scale optimisation problems, such as progressive hedging type method \citep{Munoz2015}. Although some decomposition algorithms have been proposed to tackle the computational difficulty and claimed to be capable of solving problems with short-term and long-term uncertainty \citep{Downward2020}, the algorithms were only demonstrated to solve a problem with only short-term \citep{Munoz2016} or long-term uncertainty \citep{Singh2009}. Therefore, this paper proposes an algorithm for solving such problems efficiently and demonstrates the algorithm for a long-term investment planning problem with short-term and long-term uncertainty. 

In this paper, we propose an algorithm to efficiently solve large-scale optimisation problems that exhibit a column bounded block-diagonal structure, where Subproblems (\textbf{SP}s) differ on the right-hand side and cost coefficients. Such problems can be formulated as a full Master Problem (\textbf{MP}) \eqref{eq:MP},

\begin{equation}
    \min_{\mathbf{x} \in \mathcal{X}} f(\mathbf{x})+\sum_{i \in \mathcal{I}}\pi_{i}g(x_{i}, c_{i}),
    \label{eq:MP}
\end{equation}
where $f(\mathbf{x})=\sum_{i \in \mathcal{I}}\pi_{i}c_{i}x_{i}$ and the function $g(x_{i},c_{i})$ is the optimal solution of the linear programming \textbf{SP}, 
\begin{equation}
    g(x_{i},c_{i}):=\min_{y_{i} \in \mathcal{Y}}\{c_{i}^{\top}Cy_{i}| Ay_{i}\leq Bx_{i}\}.
    \label{eq:SP}
\end{equation}
The set of decision nodes is given by $\mathcal{I}$. The $x_{i}$ are subvectors of $\mathbf{x}$. The $y_{i}$ is the decision variables of \textbf{SP} $i$ that is in the convex set $\mathcal{Y}$. The $\pi_{i}$ are non-negative constants. The coefficient matrices $A$, $B$, and $C$ are the same in \textbf{SP}s, and $x_{i}$ and $c_{i}$ are independent of short-term uncertainty. The decisions made in $f(\mathbf{x})$ are passed to \textbf{SP}s as right-hand side parameters. Multi-stage stochastic programming problems can be formulated as \eqref{eq:MP} and \eqref{eq:SP}, but the algorithm can be applied to any optimisation problems with the same structure.

Solving this problem directly can be computationally expensive. However, when $g(x_i,c_i)$ is convex and decreasing w.r.t. $x_i$, and concave and increasing w.r.t. $c_i$, one can exploit these properties to  efficiently solve the problem \citep{Mazzi2020}. \citep{Mazzi2020} proposed two inexact oracles that approximate $g(x_i,c_i)$ from below and above adaptively, and by using these, one can avoid solving all \textbf{SP}s every iteration to reduce the computational cost compared with standard Benders decomposition. The method is called Adaptive Benders. Like other Benders-type decomposition, Adaptive Benders suffers from oscillation and results in slower performance. The performance of the Adaptive Benders method was tested on a UK power system planning problem \citep{Mazzi2020}. However, it is a single-region investment planning problem, and we found that the algorithm becomes slower after introducing more regions into the problem. This issue needs to be addressed because a power system investment planning problem normally involves multiple regions connected via transmission lines \citep{Gacitua2018}. Therefore, this paper develops a stabilised Benders decomposition with adaptive oracles. We call the improved method stabilised Adaptive Benders in the rest of the paper. The stabilised Adaptive Benders consist of a level set method stabilisation and a mechanism that dynamically selects the \textbf{SP}s to solve at every iteration. 

We use the algorithm to solve an investment planning problem with short-term and long-term uncertainty formulated as multi-horizon stochastic programming. In such a problem,  $\mathbf{x}$ represents investment decisions with corresponding investment cost $f(\mathbf{x})$. The investments affect a set $\mathcal{I}$ of investment periods, and $x_{i}$ is the subvector of $\mathbf{x}$ that represents the investments that affect period $i$, $c_{i}$ specifies the operational costs, $y_{i}$ defines the operational decisions at period $i$, and $g(x_{i}, c_{i})$ gives the optimal operational cost. The $\pi_i$ is the probability associated with decision node $i$.

The contributions of this paper are: (1) we develop a level set stabilised Adaptive Benders decomposition to address the oscillation issue and analyse the tuning of parameters; and (2) we test the proposed method on a multi-horizon stochastic programming model with short-term and long-term uncertainty with up to $1$ billion variables and $4.5$ billion constraints. The results show that it is up to $113.7$ times faster than standard Benders and $2.14$ times faster than the Adaptive Benders for a $1.00\%$ convergence, and up to $45.5$ times faster than standard Benders for a $0.10\%$ convergence and the unstabilised Adaptive Benders cannot solve the largest instance to convergence tolerance; and (3) dynamic level set stabilisation increases the robustness of the proposed method and can be up to  $33.5$ times faster for 1.00\% convergence and $25.4$ times faster for 0.10\% convergence compared with standard level set method stabilisation with poor parameter choices.

The outline of the paper is as follows: Section \ref{sec:literature_review} introduces the background knowledge regarding stochastic programming, multi-horizon modelling approach, Benders decomposition and stabilisation. Section \ref{sec:LM} introduces the level set stabilisation method. Section \ref{sec:problem_description_modelling_assumptions} gives the problem description. Section \ref{sec:model} presents the model for the case study. Section \ref{sec:results} states the computational results and numerical analysis. Section \ref{sec:discussion} discusses the implications of the method and results and summaries the limitations of the research. Section \ref{sec:conclusions} concludes the paper and suggests further research.

\section{Literature review}
\label{sec:literature_review}
This paper proposes a Benders-type algorithm to solve large-scale optimisation problems. In the following, we present the background knowledge of stochastic programming, multi-horizon modelling approach, standard Benders decomposition, Adaptive Benders decomposition, and level set method stabilisation.

\subsection{Stochastic programming}
\label{sec:stochastic_programming}
Stochastic programming is the part of mathematical programming and operations research that studies how to incorporate uncertainty into decision problems \citep{King2012}. It is one of the most popular methods of dealing with uncertainties in energy system planning \citep{Birge2011}. The electricity system in regulated markets is a well-developed area for using stochastic programming in energy \citep{Wallace2003,powell2016}. However, stochastic programming is also exploited in natural gas systems \citep{Fodstad2016}, offshore oil and gas infrastructure planning \citep{Gupta2014}, and hydrogen network \citep{Galan2019}.

Two-stage stochastic programming \citep{boffino2019}, multi-stage stochastic programming \citep{pereira1991}, stochastic mixed-integer programming \citep{salo2022, lara2020,Munoz2016}, and stochastic nonlinear programming \citep{li2021phd} are all used in energy system research. In \citep{lara2020}, a multi-stage stochastic mixed-integer programming formulation was developed to optimise electricity infrastructure planning over multiple years. In order to solve a large-scale model, they decomposed and solved the problem using parallelised stochastic dual dynamic integer programming. 

\subsection{Multi-horizon stochastic programming}
\label{sec:multi_horizon_programming}
\begin{figure}[t]
    \centering
    \includegraphics[scale=1]{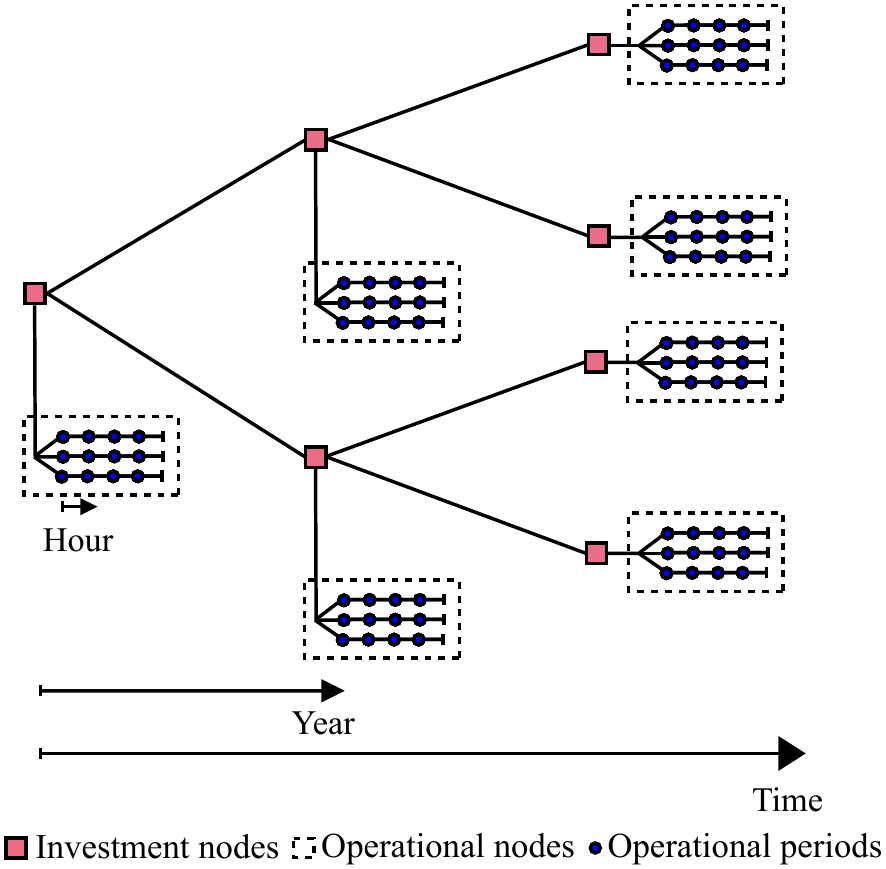}
    \caption{Illustration of multi-horizon stochastic programming with short-term and long-term uncertainty.}
    \label{fig:multi-horizon_tree}
\end{figure}
In traditional multi-stage stochastic programming, uncertainty from operational and strategic time horizons can lead to a large scenario tree, thus, an intractable planning model. The multi-horizon modelling approach was proposed as an alternative formulation that reduces the model size significantly \citep{Kaut2014}. One can have a much smaller model by disconnecting operational nodes between successive planning stages and embedding them into their respective strategic nodes. The resulted model is called multi-horizon stochastic programming. However, the multi-horizon formulation is an approximation to multi-stage stochastic programming unless two requirements are met \citep{Kaut2014}: a) strategic and operational uncertainties are independent, and the strategic decisions must not depend on any particular operational decisions; and b) the operational decisions in the last operational period in a stage do not affect the system operation in the first operational period in the next stage. An illustration of multi-horizon stochastic programming with short-term and long-term uncertainty is shown in Figure \ref{fig:multi-horizon_tree}.

\subsection{Benders decomposition}
\label{sec:benders_decomposition}
Benders decomposition was firstly developed in \citep{Benders1962} and has been successfully applied to a wide range of difficult optimisation problems \citep{Rahmaniani2017}. Benders decomposition exploits the block diagonal structure of \eqref{eq:MP} and creates outer linearisation. This method has been extended in stochastic programming to take care of feasibility questions and is known as the L-shaped method \citep{VanSlyke1969}. In Benders decomposition, a sequence of approximations is solved, and two types of constraints are added after each solve: feasibility cuts (enforcing the feasibility of \eqref{eq:MP}) and optimality cuts (linear approximations to \eqref{eq:MP} on its domain of finiteness)\citep{Birge2011}.

In standard Benders decomposition, a relaxation of the \textbf{MP} is solved. At iteration $j$, the Relaxed Master Problem (\textbf{RMP}) is
\begin{subequations}
\begin{align}
&\min_{\mathbf{x} \in \mathcal{X}, \beta} f(\mathbf{x}) + \sum_{i \in \mathcal{I}}\pi_{i}\beta_{i}\\
&\text{s.t. }\beta_{i} \geq \theta + \lambda^{\top}(x_{i}-x),\hspace{1cm}  (x,\theta,\lambda) \in F_{i(j-1)}, i \in \mathcal{I},
\end{align}
\end{subequations}
where $F_{i(j-1)}$ is the set of cuts associated with \textbf{SP} $i$ up to iteration $j-1$. To perform Benders decomposition, we firstly solve the \textbf{RMP} to obtain optimal solution $\mathbf{x}_j$. Then we pass a  subvector of $\mathbf{x}_j$, $x_{ij}$, to the \textbf{SP} $i$ and call an oracle that gives the optimal value of the \textbf{SP}, $\theta_{ij}$, and a subgradient, $\lambda_{ij}$,  w.r.t $x_{ij}$.  Finally, a new cutting plane is added to $F_{i(j-1)}$ which gives $F_{ij}:=F_{i(j-1)}\cup \{x_{ij},\theta_{ij},\lambda_{ij}\}$
The standard Benders decomposition is presented in Algorithm \ref{alg:stand_benders}.
\begin{algorithm}[t]
\caption{Standard Benders}\label{alg:stand_benders}
\begin{algorithmic}[1]
     \State choose $\epsilon$ (convergence tolerance), $\underline{\beta}$ (initial lower bound for $\beta_{i}$), $U^{*}_0:=M$ (initial upper bound), set  $j:=0$, $F_{i0}:=\{(\beta_{i0},0,0)\}$ for each $i \in \mathcal{I}$;
     \Repeat
      \State set $j:=j+1$;
      \State solve \textbf{RMP} and obtain $\beta_{ij}$ and $\mathbf{x}^{RMP}_{j}$; set $L^{*}_{j}:=f(\mathbf{x}^{RMP}_{j})+\sum_{i \in \mathcal{I}}\pi_{i}\beta_{ij}$;
      \For{$i \in \mathcal{I}$}
      \StatePar{solve \textbf{SP} $i$ at $(x^{RMP}_{ij},c_{i})$ and  obtain $\theta_{ij}$ and $\lambda_{ij}$;}
      \EndFor
     \StatePar{set $U^{*}_{j}:=\min(U^{*}_{j-1},f(\mathbf{x}^{RMP}_{j})+\sum_{i \in \mathcal{I}}\pi_{i}\theta_{ij})$;}
        \For{$i \in \mathcal{I}$}
      \StatePar{set $F_{ij}:=F_{i(j-1)}\cup \{(x^{RMP}_{ij},\theta_{ij},\lambda_{ij})$;}
      \EndFor
     \Until{$U^{*}_{j}-L^{*}_{j} \leq \epsilon.$}
    %   \algstore{testCont}
\end{algorithmic}
\end{algorithm}

\subsection{Benders decomposition with adaptive oracles}
\label{sec:benders_ao}
An investment planning problem formulated as \eqref{eq:MP} and \eqref{eq:SP} can easily get intractable once we have a large number of decision nodes that are caused by a long planning horizon or inclusion of multiple uncertainties. In stochastic programming, it refers to the curse of dimensionality that arises from the number of nodes in a scenario tree as the number of scenarios increases \citep{Powell2011}. Benders-type algorithms iteratively approximate the \textbf{SP} cost function through a set of cutting planes. However, the acquisition of the cutting planes needs all \textbf{SP}s to be solved at every iteration. Thus, Benders decomposition may get slow severely when there are many \textbf{SP}s. Therefore, research on making Benders decomposition more efficient was conducted \citep{Skar2014,Zakeri,Baena2020a}. One approach is to exploit the \textbf{SP} structure to avoid solving all \textbf{SP}s but still get a valid cutting plane each iteration.

In \citep{Mazzi2020}, two adaptive oracles were proposed in order to approximate the unsolved \textbf{SP}s objective function using the solution from the solved ones, which improves the efficiency significantly. One adaptive oracle generates inexact but valid cutting planes, and the other adaptive oracle gives a valid upper bound of the actual optimal value. 

At iteration $j$, the \textbf{RMP} is
\begin{subequations}
\begin{align}
&\min_{\mathbf{x} \in \mathcal{X}, \beta} f(\mathbf{x}) + \sum_{i \in \mathcal{I}}\pi_{i}\beta_{i}\\
&\text{s.t. }\beta_{i} \geq \underline{\theta} + \underline{\lambda}^{\top}(x_{i}-x),\hspace{1cm} (x,\underline{\theta},\underline{\lambda}) \in F_{i(j-1)}, i \in \mathcal{I},
\end{align}
\end{subequations}
where $F_{i(j-1)}$ is the supporting hyperplane from solution of \textbf{SP} $i$ at iteration $j$. At iteration $j$, the lower bound is denoted as $\mathit{L}_{j}^{*}=f(\mathbf{x})+\sum_{i \in \mathcal{I}}\beta_{i}$ and the upper bound is  $\mathit{U}_{j}^{*}=f(\mathbf{x})+\sum_{i \in \mathcal{I}}\pi_{i}\overline{\theta}_{i}$. $\underline{\theta}$, $\underline{\lambda}$ are obtained by solving the lower bound oracle and $\overline{\theta}$ is obtained by computing the upper bound oracle. Their algorithm requires \textbf{SP} to be convex and decreasing w.r.t $x_{i}$ and concave and increasing w.r.t. $c_{i}$. In order to apply the algorithm, the \textbf{SP}s need to be always feasible, which can be achieved by introducing penalty terms in \textbf{SP}s. 

\section{Level set method stabilisation}
\label{sec:LM}
In this paper, we stabilise the algorithm in \citep{Mazzi2020} using the level set method. The level set method was introduced in \citep{Lemarechal1995}. It was then used to regularise standard Benders decomposition \citep{Fabian2000}. 

\begin{figure*}[!htb]
    \centering
    \begin{subfigure}[]{0.475\textwidth}
        \centering
        \includegraphics[scale=1]{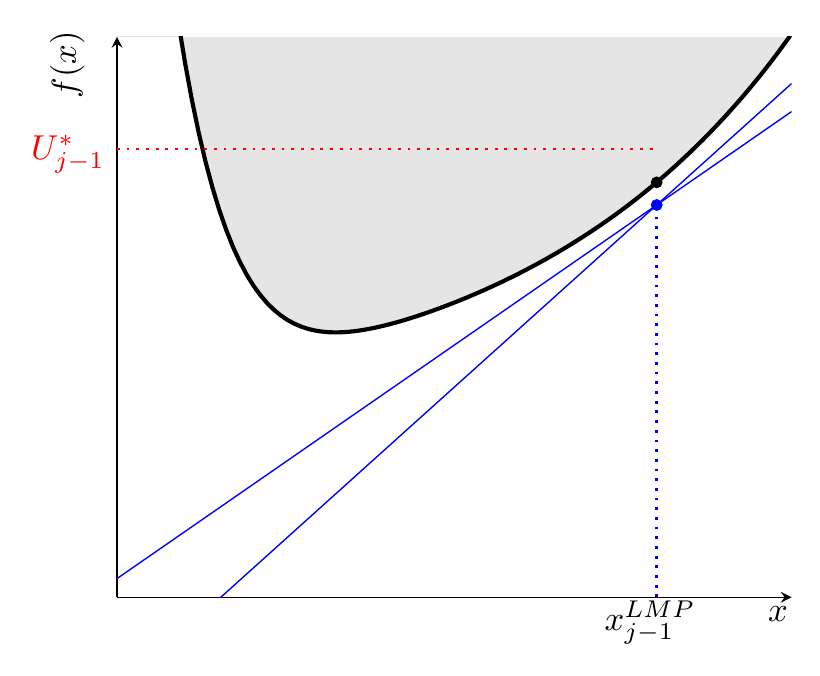}
        \caption{\small The beginning of iteration $j$}
    \end{subfigure}
    \hfill
    \begin{subfigure}[]{0.475\textwidth}
        \centering
        \includegraphics[scale=1]{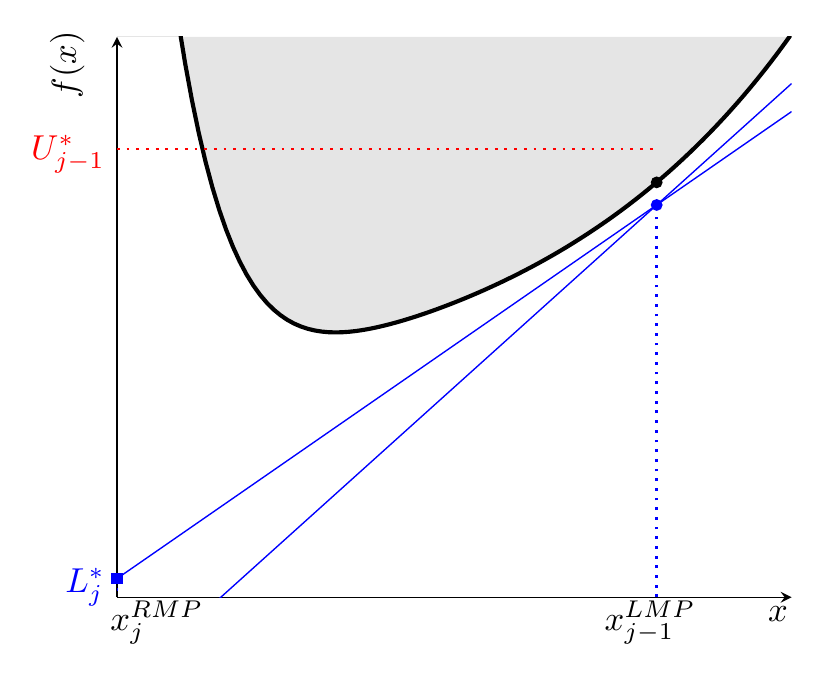}
        \caption{\small Solve the \textbf{RMP} and get a lower bound}
    \end{subfigure}
    \vskip\baselineskip
    \begin{subfigure}[]{0.475\textwidth}
        \centering
        \includegraphics[scale=1]{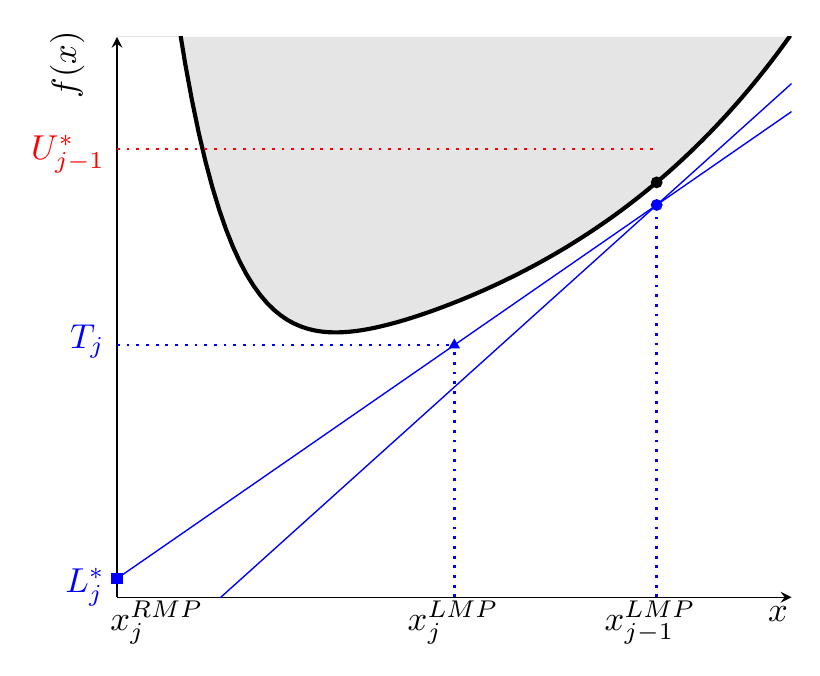}
        \caption{\small Set the target and solve the \textbf{LMP}}
    \end{subfigure}
    \hfill
    \begin{subfigure}[]{0.475\textwidth}
        \centering
        \includegraphics[scale=1]{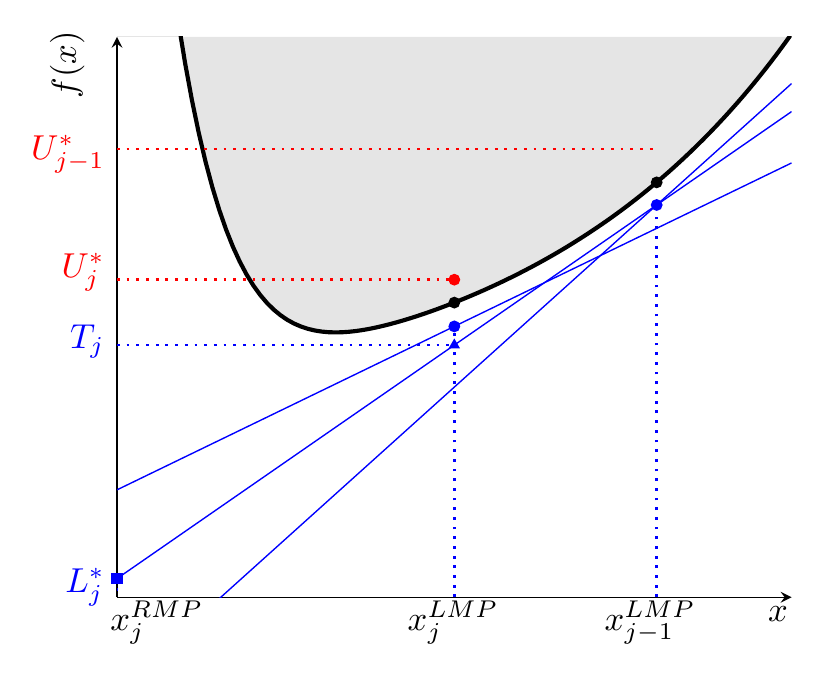}
        \caption{\small Get an upper bound and add a cut}
    \end{subfigure}
    \caption{\small An illustrative example from iteration $j-1$ to iteration $j$.}
    \label{fig:algorithm_illustration}
\end{figure*}

We now present the stabilisation step and its coordination with Adaptive Benders. At each iteration $j$, the Level Method Problem (\textbf{LMP}) for stabilisation can be formulated as
\begin{subequations}
\begin{align}
    &\min_{\mathbf{x}\in \mathcal{X},\beta} \norm{\mathbf{x}-\mathbf{x}_{j-1}}^2\\
    &\text{s.t. } \beta_{i} \geq \underline{\theta} + \underline{\lambda}^{\top}(x_{i}-x),
    \hspace{2cm}  (x,\underline{\theta},\underline{\lambda}) \in F_{i(j-1)}, i \in \mathcal{I}\\
    &\phantom{s.t. }f(\mathbf{x})+\sum_{i \in \mathcal{I}}\pi_{i}\beta_{i}\leq L^{*}_{j}+\gamma\Delta_{j}. \label{eq:level_set}
\end{align}
\end{subequations}
Constraint \eqref{eq:level_set} is the level set, $L^{*}_{j}+\gamma\Delta_{j}$ is the target that is denoted as $T_{j}$, and $\Delta_{j}=U^{*}_{j-1}-L^{*}_{j}$. The stabilisation factor, $\gamma$, is interpreted as the ratio of the achieved improvement to the predicted improvement between successive iterations. The lower bound is denoted by $\mathit{L}^{*}_{j}$ is the lower bound. By introducing \textbf{LMP}, we restrict expected improvement between iterations, thus restricting the distance moved between iterations.
% Thus we can add the target $\mathit{T}=f(\mathbf{x})+\sum_{i \in I}\beta_{i}$ to the approximated problem, where $f(\mathbf{x})$ is the investment cost based on the decision from \textbf{LMP}.
\textbf{LMP} is essentially a \textbf{RMP} but with an objective that minimises distance and an extra constraint for the level set.
A graphical interpretation of the level method stabilised Adaptive Benders decomposition is presented in Figure \ref{fig:algorithm_illustration}.

Figure \ref{fig:algorithm_illustration} illustrates how the algorithm works for iteration $j$. At the beginning of iteration $j$, we have the cuts that have been added in all previous iterations, the upper bound $U^{*}_{j-1}$ and the lower approximation of the function value (blue dot). The black dot represents the function value at that point which is unknown unless all \textbf{SP}s are solved exactly. Based on the cuts, we solve the \textbf{RMP} and get a lower bound $L^{*}_{j}$ (blue square). If there is no stabilisation, we would move to the point $x^{RMP}_{j}$. When there is stabilisation, the moving area is restricted by the target $T_{j}$, and we move to the closest point to $x^{LMP}_{j-1}$ that is below the target.  At point $x^{LMP}_{j}$, we evaluate one or more \textbf{SP}s and get a new upper bound $U^{*}_{j}$ (red dot) and add cuts. 

\begin{algorithm}[!htb]
    \caption{Level set method stabilised Benders decomposition with adaptive oracles}\label{alg:Level set method stabilised Benders decomposition with adaptive oracles}
    \begin{algorithmic}[1]
        \State choose $\epsilon$ (convergence tolerance), $\gamma$ (stabilisation factor), $\underline{\beta}$ (initial lower bound $\beta_{i}$), $U^{*}_0:=M$ (initial upper bound);
        \StatePar{set  $j:=0$, $F_{i0}:=\{(\beta_{i0},0,0)\}$ for each $i \in \mathcal{I}$;}
        \State solve \textbf{SP} at the special point $(\underline{x},\underline{c})$ and obtain $\theta$, $\lambda$ and $\phi$; set $\mathcal{S}:=\{(\underline{x},\underline{c},\theta,\lambda,\phi)\}$;
        \Repeat
            \State set $j:=j+1$;
            \State solve \textbf{RMP} and obtain $\beta_{ij}$ and $\mathbf{x}^{RMP}_{j}$; set $L^{*}_{j}:=f(\mathbf{x}^{RMP}_{j})+\sum_{i \in \mathcal{I}}\pi_{i}\beta_{ij}$;
            \StatePar{set $\mathbf{x}^{Ref}:=\mathbf{x}^{RMP}_{j}$ (when $j=1$); }
            \StatePar{set \textbf{LMP} target: $L^{*}_{j}+\gamma(U^{*}_{j-1}-L^{*}_{j})$;}
            \StatePar{solve \textbf{LMP} and obtain $\mathbf{x}^{LMP}_{j}$;}
            \For{$i \in \mathcal{I}$}
                \StatePar{call adaptive oracles at $(x^{LMP}_{ij},c_{i})$ and  obtain $\underline{\theta}_{ij}$, $\overline{\theta}_{ij}$, $\overline{\phi}_{ij}$ and $\underline{\lambda}_{ij}$;}
            \EndFor
            \State set $n:=0$;
            \Repeat
                \State $n:=n+1;$\;
                \State choose $\textbf{SP}_i$ in $\mathcal{I}$ that has the largest gap;
                \State solve $\textbf{SP}_i$ at $(x^{LMP}_{ij},c_{i})$ exactly and obtain $\theta_{ij}$, $\lambda_{ij}$, $\phi_{ij}$;
                \State set $\mathcal{S}:=\mathcal{S}\cup\{(x^{LMP}_{ij},c_{i},\theta_{ij},\lambda_{ij},\phi_{ij})\}$;
                \For{$i \in \mathcal{I}$}
                    \StatePar{set $F_{ij}:=F_{i(j-1)}\cup \{(x^{LMP}_{ij},\underline{\theta}_{ij},\underline{\lambda}_{ij})$;}
                \EndFor
                \For{$i \in \mathcal{I}$}
                    \StatePar{call adaptive oracles at $(x^{LMP}_{ij},c_{i})$ and obtain $\underline{\theta}_{ij}$, $\overline{\theta}_{ij}$, $\overline{\phi}_{ij}$ and $\underline{\lambda}_{ij}$;}
                \EndFor
                \StatePar{set $L^{LBO}_{j}:= f(\mathbf{x}^{LMP}_{j}) + \sum_{i \in \mathcal{I}}\pi_{i}\underline{\theta}_{ij}$;}
                \StatePar{set $U^{UBO}_{j}:= f(\mathbf{x}^{LMP}_{j}) + \sum_{i \in \mathcal{I}}\pi_{i}\overline{\theta}_{ij}$;}
            \Until{$U^{UBO}_{j}-L^{LBO}_{j} \leq U^{*}_{j-1}-L^{*}_{j-1} $ \normalfont{\textbf{or}} $n>|\mathcal{I}|$ \normalfont{\textbf{or}} $L^{LBO}_{j}\geq U^{*}_{j-1};$} 
            \StatePar{set $U^{*}_{j}:=\min(U^{*}_{j-1},U^{UBO}_{j})$, $\mathbf{x}^{Ref}:=\mathbf{x}^{LMP}_{j};$}
        \Until{$U^{*}_{j}-L^{*}_{j} \leq \epsilon.$}
    %   \algstore{testCont}
    \end{algorithmic}
\end{algorithm}
\subsection{Stabilised Benders decomposition with adaptive oracles algorithm}
\label{sec:algorithm}
In this section, we present the stabilised Benders decomposition with the adaptive oracles algorithm, shown in Algorithm \ref{alg:Level set method stabilised Benders decomposition with adaptive oracles}. 

Although an interesting feature of the level method compared with other bundle-type methods is that $\gamma$ is fixed \citep{Zverovich2012}, it may be beneficial to adjust it based on the progress. Therefore, in addition to using fixed stabilisation factor $\gamma$, we also explore adjusting the stabilisation factor for a potentially more robust algorithm. There are several ways to adjust stabilisation dynamically, and we choose a method analogous to what is used to adjust trust regions. The trust region method uses a local approximation of the function to be minimised and optimised within the trust region. The trust region size is updated throughout the iterations. In the trust region method, one adjusts the trust region according to the ratio of the actual decrease to the predicted decrease \citep{fletcher2000}. Inspired by \citep{fletcher2000}, we adjust the level set based on the ratio of the actual improvement to the expected improvement. At iteration $j$, we define the ratio
\begin{equation}
    \label{ratio}
    r:=\frac{L^{LBO}_{j-1}-L^{LBO}_{j}}{L^{LBO}_{j-1}-{T_{j}}},
\end{equation}
where $I^{A}_{j}=L^{LBO}_{j-1}-L^{LBO}_{j}$ is the actual improvement from iteration $j-1$ to $j$, and $I^{P}_{j}=L^{LBO}_{j-1}-{T_{j}}$ is the predicted improvement from iteration $j-1$ to $j$. Then we update $\gamma$, \\
\begin{algorithmic}[!htb]
    % \State choose $\omega$, $\underline{P}$ , $\overline{P}$  set;
    \If{$I^{A}_{j}>0$ \textbf{and} $I^{P}_{j}>0$}{}
        \If{$r \leq \underline{P}$}{}
        \State{$\gamma=1-\omega(1-\gamma)$}
        \ElsIf{$\underline{P} <r< \overline{P}$}{}
        \State{$\gamma=\gamma$}
        \Else
        \State{$\gamma=\omega\gamma$}
        \EndIf
    \Else{$I^{A}_{j}\leq0$ \textbf{or} $I^{P}_{j}\leq0$}{}
        \State{Inexact information, $\gamma=\gamma$}
    \EndIf
\end{algorithmic}
where $0 \leq \omega \leq 1$ is a constant that increases or decreases $\gamma$, $\overline{P}$ and $\underline{P}$ are constants that determine what actions to take on $\gamma$. Unlike standard Benders that knows the exact value of the \textbf{SP}s, stabilised Adaptive Benders only knows lower and upper bounds on the objective values. By comparing the lower and upper bounds with the exact values of the \textbf{SP}s, we find that the lower bound oracle gives a much closer and more stable approximation. Therefore, we use $L^{LBO}$ instead of $U^{UBO}$ when defining the ratio $r$. Furthermore, a bad approximation from the lower bound oracle at the current point or a bad approximation from the upper bound oracle from previous points can lead to a negative $I^{P}_{j}$. Therefore, we choose to do nothing once we find that $I^{P}_{j}$ is negative. For $I^{A}_{j}$, there are two possibilities for it to be negative: 1) bad approximation from the lower bound oracle or the upper bound oracle, and 2) going to a bad point. A bad point, in this case, means the $L^{LBO}$ at current iteration is higher than the $L^{LBO}$ at the previous iteration. If the information is exact and $I^{A}_{j}$ is negative, one may reject the point, go back to the best point seen so far, and try again with a higher $\gamma$. However, in the case of inexact information, it may not be sensible to reject a point based on a bad approximation. In the computational study in this paper, fixed stabilisation is mainly used, but dynamic stabilisation is also tested. 

\section{Problem description, modelling strategies and modelling assumptions}
\label{sec:problem_description_modelling_assumptions}
The proposed power planning problem is designed to choose the optimal investment strategy and operating scheduling for a power system to achieve emission targets. In this section, we present the temporal and geographical representations of the problem and the modelling assumptions. 

The problem under consideration aims to make optimal investment and operational decisions for the UK power system that satisfies the emission reduction goal under a) short-term uncertainty, including renewable energy availability and load profile; and b) long-term uncertainty, including CO$_2$ budget, CO$_2$ tax, and long-term power demand. 

For the investment planning, we consider: (a) thermal generators (Coal-fired plant, OCGT, CCGT, Diesel, and nuclear plants); (b) generators with Carbon Capture and Storage (CCS) (Coal-fired plant with CCS); (c) renewable generators (offshore wind, onshore wind and solar PV); (d) electric storage (PHES and lithium); and (e) transmission lines. The capital expenditures and fixed operational costs are assumed to be known. The problem is to determine: (a) the capacities of technologies and (b) operational strategies that include scheduling of generators, storage and approximate power flow among regions to meet the power demand with minimum overall investment, operational and environmental costs.

\subsection{Modelling strategies and assumptions}
In this section, we present the modelling strategies and assumptions we use in the stochastic long-term multi-region multi-period investment planning problem.

\subsubsection{Scenario generation}
For short-term uncertainty, we select some time intervals with a half-hourly resolution in four seasons of a year and scale them up to represent an operational scenario. For long-term uncertainty, each independent uncertain parameter has $n$ possible outcomes in the next stage, which is linked to additional $n$ possible outcomes in the following stage. The realisations in one stage are assigned with an equal probability. We use a reasonably simple scenario generation routine because scenario generation is not the scope of the paper, and we refer the readers to \citep{King2012,Fairbrother2022} for more advanced scenario generation approaches. 

\subsubsection{Geographical representation of the problem}
The problem potentially consists of many regions and results in a large model. Therefore, we aggregate regions into representative ones to reduce the number of locations. The generators and storage units in one region with the same characteristics are aggregated into clusters. In such a way, the model does not invest in a specific unit but in that type of device, and a linear investment model may be sufficient in this case.

\subsubsection{Modelling assumptions}
We assume that: a) a linear cost model for each technology because we deal with an aggregated system and the fixed part of the investment cost can be evened out and lead to a linear programming master problem; b) the Kirchhoff voltage law is omitted, and c) no loss in the transmission lines.  

\section{Mathematical model}
\label{sec:model}
This section presents the mathematical model for the power system investment planning and operational problem. The problem is decomposed by having an investment planning master problem and an operational \textbf{SP}. The complete nomenclature of the model can be found in \ref{nomenclature}. We use the conventions that calligraphic capitalised Roman letters denote sets, upper case Roman and lower case Greek letters denote parameters, and lower case Roman letters denote variables. The indices are subscripts, and name extensions are superscripts. The same lead symbol represents the same type of thing. The names of variables, parameters, sets and indices are single symbols. 

\subsection{Investment planning model}
\begin{subequations}
\label{mod:investment_planning}
\begin{alignat}{3}
    &\min \mathrlap{c^{INV}+\kappa\sum_{i \in \mathcal{I}}\pi_{i}c^{OPE}(x_i,c_i)} & \label{objective}\\
    &\text{s.t.}&& \mathrlap{c^{INV}=
    \sum_{i\in \mathcal{I}_{0}}\delta^{I_{0}}_{i}\pi^{I_{0}}_{i}\sum_{p \in \mathcal{P}} C^{Inv}_{pi}x^{Inst}_{pi}+\kappa \sum_{i \in \mathcal{I}}\delta^{I}_{i}\pi^{I}_{i}\sum_{p \in \mathcal{P}}C^{Fix}_{pi}x^{Acc}_{pi}} & \label{investment_cost}\\
    &&&x^{Acc}_{pi}=X^{Hist}_{p}+\sum_{i \in\mathcal{I}_{i}|\kappa(i-i_{0})\leq H^{P}_{p}}x_{pi}^{Inst},\phantom{abcdefghijklmnop}   & p \in \mathcal{P},  i \in \mathcal{I}\label{cap_tech}\\
    &&&x^{Acc}_{pi} \leq X^{Max}_{p},  & p \in \mathcal{P},  i \in \mathcal{I} \label{max_tech}\\
    &&&x_i=\left( \{x^{Acc}_{pi}, p \in \mathcal{P}\}, \mu^{DP}_{i}, \mu^{E}_{i}\right),  & i \in \mathcal{I} \label{rhs}\\
    &&&c_i=\left(C^{CO_2}_{i}\right),  & i \in \mathcal{I} \label{cost_coefficient}\\
    &&&\mathrlap{x^{Inst}_{pi}, x^{Acc}_{pi} \in \mathbb{R}^{+}_{0}.} \label{domain1}
\end{alignat}
\end{subequations}

The total cost for investment planning, Equation \eqref{objective}, consists of actual discounted investment costs and discounted fixed operating and maintenance costs $c^{INV}$, as well as the expected operational cost of the system over the time horizon $\kappa\sum_{i \in \mathcal{I}}\pi_{i}c^{OPE}(x_i,c_i)$. Here, $\kappa$ is a scaling factor that depends on the time step between two successive investment nodes. Constraint \eqref{cap_tech} states that the accumulated capacity of a technology $x^{Acc}_{pi}$ in an operational node equals the sum of the historical capacity $X^{Hist}_{p}$ and newly invested capacities $x^{Inst}_{pi}$ in its ancestor investment nodes $\mathcal{I}_{i}$ that are in their lifetimes. The parameter $X^{Max}_{p}$ denotes the maximum accumulated capacity of technologies. Constraint \eqref{rhs} collects all right hand side coefficients that will be passed to the \textbf{SP} \eqref{mod:operational} into vector $x_i$. And constraint \eqref{cost_coefficient} collects all the cost coefficients into vector $c_i$.

\subsection{Operational model}
We now compute the operational cost $c^{OPE}(x_i,c_i)$ at one operational node $i$ by solving \textbf{SP} \eqref{mod:operational} given the decisions $x_i$ and $c_i$ made in the master problem \eqref{mod:investment_planning}. Note that we omit index $i$ in the operational model for ease of notation.
\begin{subequations}
\label{mod:operational}
\begin{alignat}{3}
    &\min \quad \mathrlap{\sum_{t \in \mathcal{T}}\pi_{t}H_{t}\left(\sum_{g \in \mathcal{G}}C^{G}_{g}p^{G}_{gt}+ \sum_{s \in \mathcal{S}}C^{S}_{s}p^{SE+}_{st}\sum_{z \in \mathcal{Z}}C^{Shed}p^{ShedP}_{zt} \right)}& \label{SP_objective}\\
    &\text{s.t.} && p^{G}_{gt} \leq p^{AccG}_{g}, & g \in \mathcal{G}, t \in \mathcal{T}\label{gmod:tech_cap}\\
    &&& -p^{AccL}_{l} \leq p^{L}_{lt} \leq p^{AccL}_{l}, & l \in \mathcal{L}, t \in \mathcal{T}\label{gmod:line_cap}\\
    &\quad&& p_{st}^{SE+} \leq p_{s}^{AccSE}, & s \in \mathcal{S}, t \in \mathcal{T} \label{gmod:estore_charge_cap}\\
    &\quad&&p_{st}^{SE-} \leq p_{s}^{AccSE}, & s \in \mathcal{S}, t \in \mathcal{T} \label{gmod:estore_discharge_cap}\\
    &\quad&&q_{st}^{SE} \leq \gamma^{SE}_{s}p_{s}^{AccSE}, & s \in \mathcal{S}, t \in \mathcal{T} \label{gmod:estore_energy_cap}\\
    &\quad&& -\alpha^{G}_{g}p_{g}^{AccG} \leq p_{gt}^{G}-p_{g(t-1)}^{G} \leq \alpha^{G}_{g}p_{g}^{AccG}, &\phantom{ab} g \in \mathcal{G}, n \in \mathcal{N}, t \in \mathcal{T}_{n} \label{gmod:turbine_ramp}\\
    &\quad&&\mathrlap{\sum_{g \in \mathcal{G}_{z}}p_{gt}^{G}+\sum_{l \in \mathcal{L}^{In}_{z}}p_{lt}^{L}+\sum_{s \in \mathcal{S}_{z}}p_{st}^{SE-}+\sum_{r \in \mathcal{R}_{z}}R^{R}_{zt}p_{r}^{AccR}+ p_{zt}^{ShedP}=}&\notag\\
    &\quad&&\phantom{abcdef}\mu^{DP}P^{DP}_{zt}+\sum_{l \in \mathcal{L}^{Out}_{z}}p_{lt}^{L}+\sum_{s \in \mathcal{S}_{z}}p_{st}^{SE+}+p_{zt}^{GShedP},& z \in \mathcal{Z}, t \in \mathcal{T}\label{gmod:kcl}\\
    &\quad&& q_{s(t+1)}^{SE}=q_{st}^{SE}+H_{t}(\eta_{s}^{SE}p_{st}^{SE+}-p_{st}^{SE-}),\phantom{abcde} & s \in \mathcal{S}, n \in \mathcal{N}, t \in \mathcal{T}_{n} \label{gmod:estore_balance}\\
    &\quad&& \sum_{t \in \mathcal{T}}\sum_{g \in \mathcal{G}} \pi_{t}H_{t}E^{G}_{g}p^{G}_{gt}\leq \mu^{E},&\label{gmod:co2budget}\\
    &&&\mathrlap{p^{L}_{lt} \in \mathbb{R}_{0},}\\
    &&&\mathrlap{p^{G}_{gt}, p^{AccG}_{g}, p^{ShedP}_{zt}, p^{SE+}_{st}, p^{SE-}_{st}, p^{AccSE}_{s}, q^{SE}_{st},p^{AccR}_{r}, p^{GShedP}_{zt} \in \mathbb{R}_{0}^{+}.} \label{gmod:domain}
\end{alignat}
\end{subequations}
The operational cost function $c^{OPE}(x, c)$ includes total operating costs of all generators and storage facilities $C^{G}_{g}p^{G}_{gt}+C^{S}_{s}p^{SE+}_{st}$ and load shedding costs $C^{Shed}p^{Shed}$. The parameters $C^{G}_{g}$ and $C^{S}_{s}$ include the variable operational cost of generators and storage. For thermal generators, $C^{G}_{g}$ also includes the fuel cost and the CO$_2$ tax charged on the emissions of generators. Constraint \eqref{gmod:tech_cap} ensures that generators are within their capacity limits. Constraint \eqref{gmod:line_cap} shows that the power flow $p^{L}_{t}$ is within the transmission capacity $p^{AccL}_{l}$. Constraints \eqref{gmod:estore_charge_cap} and \eqref{gmod:estore_discharge_cap} dictate that the power charged $p^{SE+}_{st}$ and the discharging power $p^{SE-}_{st}$ of a storage facility should be within the capacity, respectively. Constraint \eqref{gmod:estore_energy_cap} limits the energy storage level $q^{SE}_{st}$ to be within the capacity $q^{AccSE}_{s}$. Constraint \eqref{gmod:turbine_ramp} captures how fast thermal generators can ramp up or ramp down their power output, respectively. The parameters $\alpha^{G}_{g}$ is the maximum ramp rate of thermal generators. The power nodal balance, Constraint \eqref{gmod:kcl}, ensures that in one operational period $t$, the sum of total power generation of thermal generators $p^{G}_{gt}$, power discharged from all the electricity storage $p^{SE-}_{st}$, renewable generation $R^{R}_{zt}p^{AccR}_{rt}$, power transmitted to this region, and load shed $p^{ShedP}_{zt}$ equals the sum of  power demand $P^{DP}_{zt}$ power transmitted to other regions, and power generation shed $p^{GShedP}_{zt}$. The parameter $R^{R}_{zt}$ is the capacity factor of a renewable unit that is a fraction of the nameplate capacity $p^{AccR}$. The subset of a technology in region $z$ is represented by $R_{z}:=\{r \in \mathcal{R}: r \text{ is available in region } z\}$, where $\mathcal{R}$ can be replaced by other sets of technologies. Constraint \eqref{gmod:estore_balance} states that the state of charge $q^{SE}_{st}$ in period $t+1$ depends on the previous state of charge $q^{SE}_{st}$, the charged power $p^{SE+}_{st}$ and discharged power $p^{SE-}_{st}$. The parameter $\eta^{SE}_{s}$ represent the charging efficiency. Constraint \eqref{gmod:co2budget} restricts the total emission. The parameter $H_t$ is the length of the period $t$.  The parameter $\mu^{E}$ is the CO$_2$ budget. The symbol $E^{G}_{g}$ is the emission factor per unit of power generated. The capacities $p^{AccG}_{g}$, $p^{AccL}_{l}$, $p^{AccSE}_{s}$, scaling factor of demand $\mu^{DP}$ and CO$_2$ budget $\mu^E$ are passed from the master problem \eqref{mod:investment_planning} via vector $x_i$ and CO$_2$ tax that is included in cost coefficient $C^{G}_g$ is passed from master problem \eqref{mod:investment_planning} via vector $c_i$.

\section{Results}
This section firstly uses small illustrative cases to show how stabilisation helps solve multi-region investment planning problems. Then we demonstrate the proposed algorithm on larger instances and present the computational results.
\label{sec:results}
\subsection{Illustrative cases}
We use three cases to show the value of stabilisation in a multi-region investment planning problem. A summary of the four cases is presented in Table \ref{table:Summary of the illustrate cases.}. To simplify the visualisation of the results, we consider only two types of generation, OCGT and Diesel, and a one-time investment planning problem is solved. In this case, there is no difference between Adaptive Benders and standard Benders because there is only one \textbf{SP}. 

\begin{table}[!htb]
    \centering
    \caption{Summary of the illustrative cases.}
    \label{table:Summary of the illustrate cases.}
        \begin{tabular}
            {c|l}
            \toprule
            \multicolumn{1}{c}{} & \multicolumn{1}{c}{Description} \\ \hline
            Case A   & single region, two technologies to invest (OCGT and Diesel) \\
            Case B   & two unconnected regions with sizes 60\% and 40\% of case A \\
            Case C   & Case B with a transmission line with 0 initial capacity \\
            \bottomrule
        \end{tabular}
\end{table}

\begin{figure*}[!htb]
    \centering
    \begin{subfigure}[t]{0.47\textwidth}
        \centering
        \includegraphics[scale=1]{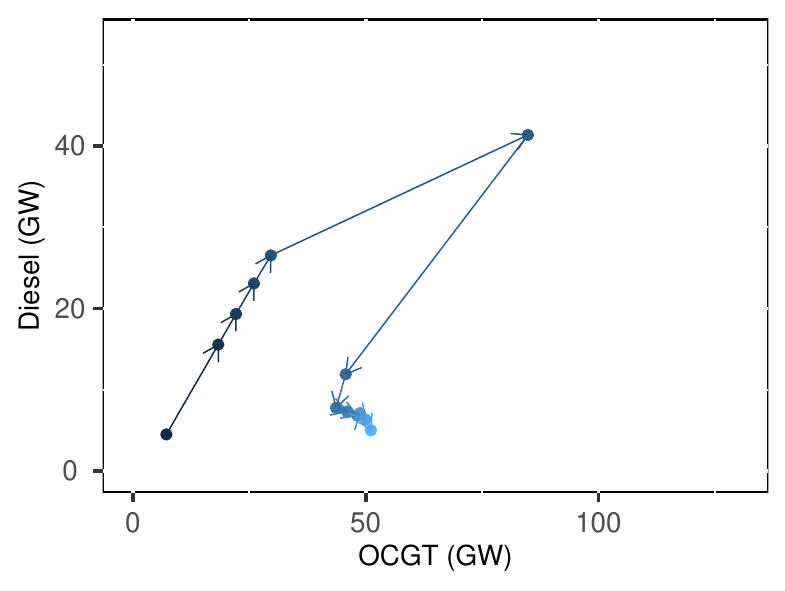}
        \caption{\small Case A (stabilised Benders, 13 iters)}
    \end{subfigure}
    \begin{subfigure}[t]{0.47\textwidth}
        \centering
        \includegraphics[scale=1]{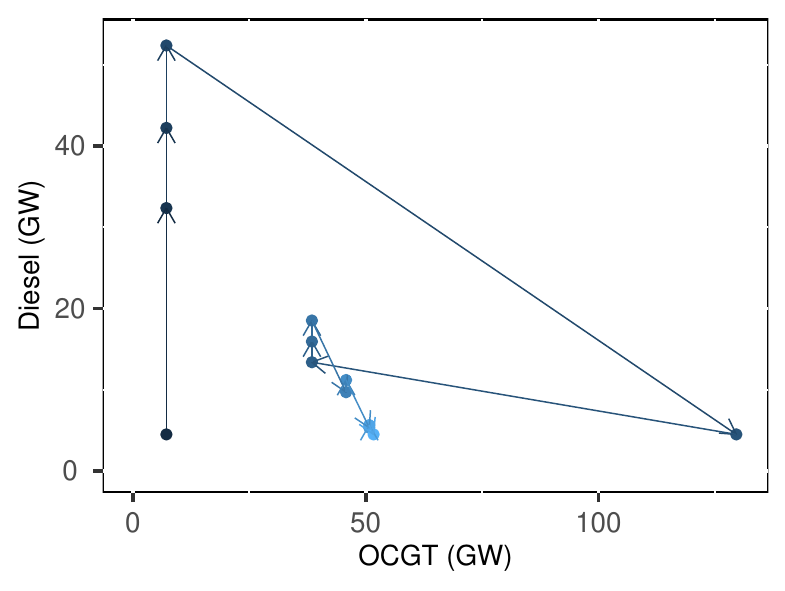}
        \caption{\small Case A (unstabilised Benders, 13 iters)}
    \end{subfigure}
    \caption{Comparative results of Case A.}
    \label{fig:Comparative results of Case A}
\end{figure*}

\begin{figure*}[!htb]
    \centering
    \begin{subfigure}[t]{0.47\textwidth}
        \centering
        \includegraphics[scale=1]{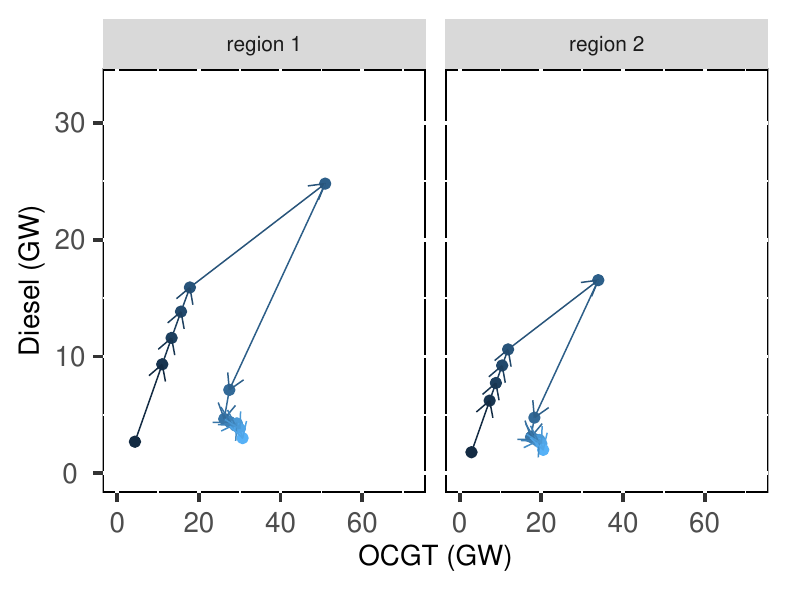}
        \caption{\small Case B (stabilised Benders, 13 iters)}
    \end{subfigure}
    \begin{subfigure}[t]{0.47\textwidth}
        \centering
        \includegraphics[scale=1]{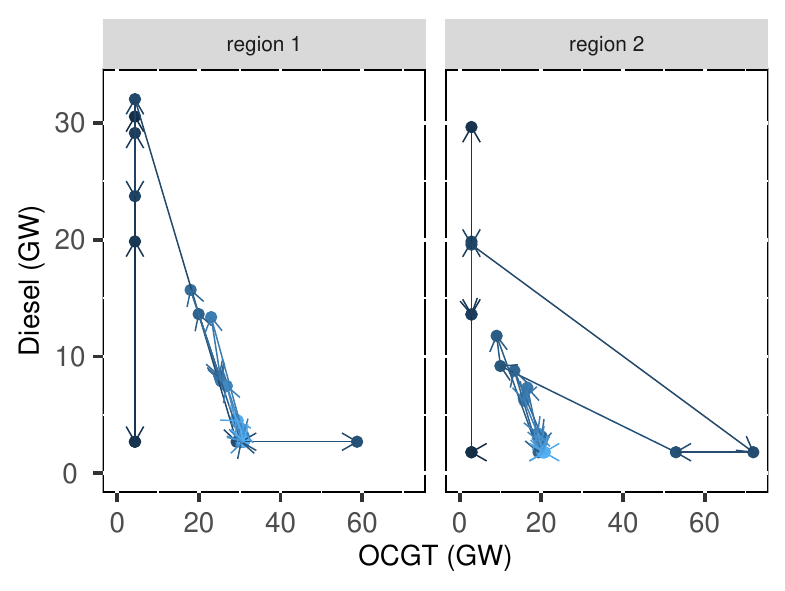}
        \caption{\small Case B (unstabilised Benders, 26 iters)}
    \end{subfigure}
    \caption{Comparative results of Case B.}
    \label{fig:Comparative results of Case B}
\end{figure*}

\begin{figure*}[!htb]
    \centering
    \begin{subfigure}[t]{0.47\textwidth}
        \centering
        \includegraphics[scale=1]{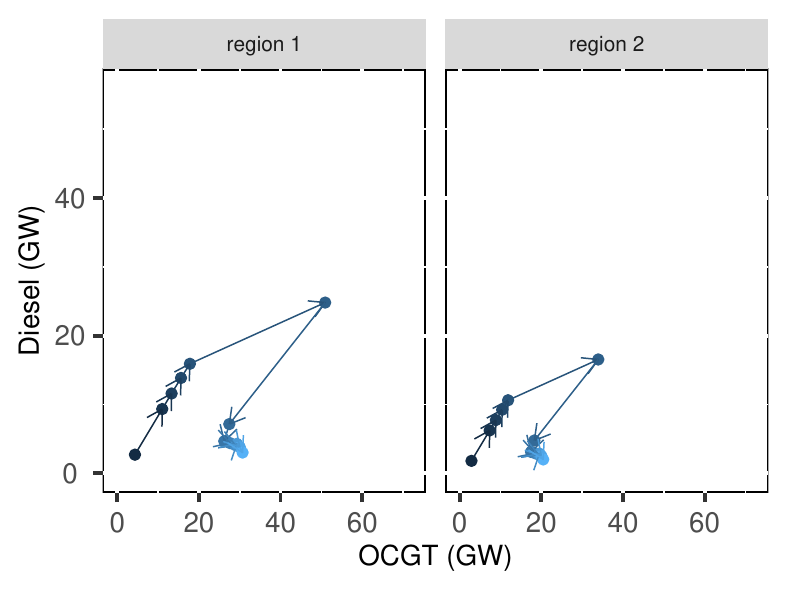}
        \caption{\small Case C (stabilised Benders, 13 iters)}
    \end{subfigure}
    \begin{subfigure}[t]{0.47\textwidth}
        \centering
        \includegraphics[scale=1]{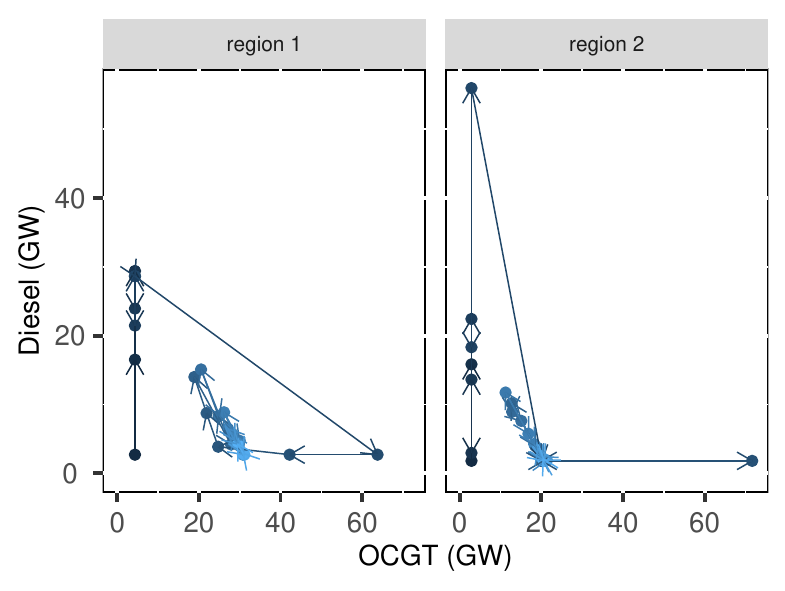}
        \caption{\small Case C (unstabilised Benders, 27 iters)}
    \end{subfigure}
    \caption{Comparative results of Case C.}
    \label{fig:Comparative results of Case C}
\end{figure*}

Figure \ref{fig:Comparative results of Case A} - Figure \ref{fig:Comparative results of Case C} show how solutions are explored until convergence. In each figure, the darkest blue point represents the initial solution, the lightest blue point is the optimal solution, and the arrows indicate the order of points explored. For the stabilised versions, the stabilisation factor is fixed to 0.2. In all cases, there is degeneracy in the dimension of the total amounts of the two generation types. From Case B and Case C, we find that there is degeneracy in the dimension of regions. In the two region cases, there is a CO$_2$ constraint that restricts the total emissions from both regions. In Case B, where to put the capacities becomes relevant. From Figure \ref{fig:Comparative results of Case B}, we see that without stabilisation, the algorithm struggles to balance the capacities of the two technologies and starts jumping to points with different proportions of the two technologies many times until it finds the optimal solution. In Figure \ref{fig:Comparative results of Case B}, we see that the stabilised approach is clearer about which direction to explore and make small movements towards the optimal instead of sampling points wildly. The number of iterations is doubled without stabilisation. For a more realistic problem with more technologies and regions and a more complicated network topology, the value of stabilisation reveals further, as is shown in Section \ref{sec:computational_results}.

In Case C, two regions are initially disconnected, but a line can be invested to connect them. However, there should be no line invested because the two regions are proportional to each other and making investments in the local generation is optimal. By observing the solution proposed by \textbf{RMP} in the unstabilised version, we notice that \textbf{RMP} does not realise that and makes an investment in the line in some iterations before finding the optimum, and this leads to more iterations compared with the stabilised version. 

\subsection{Case study}
\begin{figure}[t]
    \centering 
    \includegraphics[scale=0.8]{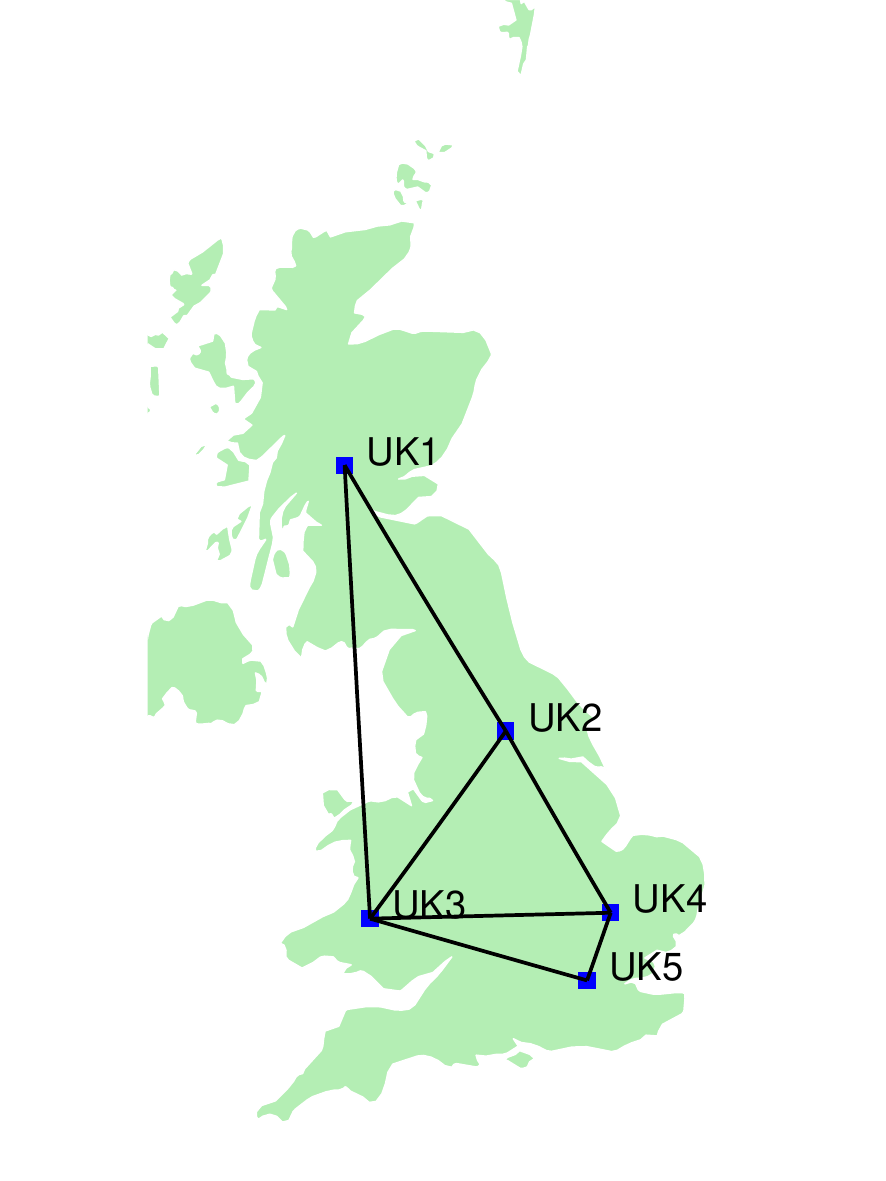}
    \caption{Illustration of the UK power system. (UK1: Scotland, UK2: North England, UK3: Midland and Wales, UK4: East England and UK5: South England.)}
    \label{fig:uk_grid}
\end{figure}

We test the stabilised Benders algorithm with adaptive oracles on the stochastic investment planning of the UK power system. We use the model presented in Section \ref{sec:model} to investigate the computational issues. The network topology is shown in Figure \ref{fig:uk_grid}. We implemented the algorithm and model in Julia 1.7.3 using JuMP \citep{jump} and solved with Gurobi 9.5.1 \citep{gurobi}.  We ran the code on nodes of a computer cluster with a 2x 3.6GHz 8 core Intel Xeon Gold 6244 CPU and 384 GB of RAM, running on CentOS Linux 7.9.2009. Some data was taken from \citep{Mazzi2020}. The Julia code and data for the case study have been made publicly available at \citep{Zhang2022_code}.

\subsection{Computational results}
\label{sec:computational_results}
This section presents the computational results of the proposed stabilised Adaptive Benders. We compare the performance of stabilised Adaptive Benders against the unstabilised versions of Adaptive Benders with one \textbf{SP} solved per iteration \citep{Mazzi2020} and the standard Benders. We use the model presented in Section \ref{sec:model} to solve a 5-region UK power system planning to make the benchmark. The long-term uncertainties include CO$_2$ price and power demand. The short-term uncertainties include wind and solar capacity factors and load profiles. The summary of cases and their problem sizes are shown in Table \ref{table:Overview of the cases used in the computational study.}. In Case 0--3, there are four short-term operational scenarios, each consisting of 4380 operational periods. Case 0 has no long-term uncertainty. Case 1 has one long-term uncertainty, CO$_2$ budget. Case 2 has CO$_2$ budget and long-term demand uncertainty. And Case 3 has CO$_2$ budget, long-term demand, and CO$_2$ tax as long-term uncertainty.

\begin{table}[t]
    \caption{Overview of the cases used in the computational study.}
%   \vspace{-0.5cm}
    \label{table:Overview of the cases used in the computational study.}
    \resizebox{\columnwidth}{!}{%
        \begin{tabular}
            {cS[table-format=3.0]S[table-format=2.0]S[table-format=2.0]|S[table-format=1.0]S[table-format=1.0]S[table-format=2.0]S[table-format=2.0]|S[table-format=1.1e1]S[table-format=1.1e1]S[table-format=1.1e1]}
            \toprule
            & {Operational periods} & \multicolumn{1}{c}{Short-term} & \multicolumn{1}{c}{Long-term} & \multicolumn{4}{c}{Number of decision nodes} & \multicolumn{3}{c}{Problem size (undecomposed)}\\
            & {per short-term scenario} & \multicolumn{1}{c}{scenarios} & \multicolumn{1}{c|}{scenarios} & \multicolumn{1}{c}{Present} & \multicolumn{1}{c}{In 5 years} & \multicolumn{1}{c}{In 10 years} & \multicolumn{1}{c|}{Total} & \multicolumn{1}{c}{Variables} & \multicolumn{1}{c}{Constraints} & \multicolumn{1}{c}{Nonzeros} \\ \hline
            Case 0  & 4380  &  4 &  1  & 1 & 1  &  1  &  3  & 2.7e6 & 7.7e6 & 1.9e7 \\
            Case 1  & 4380  &  4 &  9  & 1 & 3  &  9  & 13  & 1.6e7 & 4.6e7 & 1.1e8 \\
            Case 2  & 4380  &  4 & 81  & 1 & 9  & 81  & 91  & 1.2e8 & 3.5e8 & 8.4e8 \\
            Case 3  & 4380  &  4 & 729 & 1 & 27 & 729 & 757 &1.0e9  & 4.5e9 & 1.1e10\phantom{a}* \\
            \bottomrule
            \multicolumn{11}{r}{*: the model cannot be loaded into the solver due to its size.}
        \end{tabular}
    }
\end{table}

\begin{table}[t]
    \centering
    \caption{\centering Comparative results for standard Benders, Adaptive Benders and stabilised Adaptive Benders, $\gamma$ is fixed to 0.025 for stabilised Adaptive Benders. (speed up: the time spent using standard Benders divided by time spent using Adaptive or stabilised Adaptive Benders).}
    \label{table:Comparative results for standard Benders, adaptive Benders and stabilised adaptive Benders.}
    % \vspace{-0.45cm}
    \resizebox{\columnwidth}{!}{%
        \begin{tabular}
            {cS[table-format=1.2]|S[table-format=2.0]|S[table-format=2.0,table-column-width=3cm]S[table-format=2.0]|S[table-format=2.0,table-column-width=3cm]S[table-format=2.0]S[table-format=2.0]|S[table-format=2.0]S[table-format=5.0]S[table-format=2.1]}
            \toprule
            & \multicolumn{1}{c}{\multirow{1}{*}{$\epsilon$ (\%)}} & \multicolumn{1}{c}{Undecomposed} & \multicolumn{2}{c}{Standard Benders} & \multicolumn{3}{c}{Adaptive Benders} & \multicolumn{3}{c}{Stabilised Adaptive Benders}   \\ \cline{3-11}
            & \multicolumn{1}{c}{}  & \multicolumn{1}{c}{Time (s)} & \multicolumn{1}{c}{Iters/Evals} & \multicolumn{1}{c}{Time (s)} & \multicolumn{1}{c}{Iters/Evals} & \multicolumn{1}{c}{Time (s)} & \multicolumn{1}{c}{Speed up} & \multicolumn{1}{c}{Iters/Evals} & \multicolumn{1}{c}{Time (s)}  & {Speed up} \\
            \hline
 Case 0	&	1.00	&	440	&	18/36	&	1051	&	30/31		&	874		&	1.2		&	16/26	&	751	&	1.4	\\
	&	0.10	&		&	33/66	&	1925	&	66/67		&	1925		&	1.0		&	35/47	&	1344	&	1.4	\\
Case 1	&	1.00	&	$\infty$	&	16/192	&	5698	&	32/33		&	953		&	6.0		&	20/27	&	791	&	7.2	\\
	&	0.10	&		&	28/336	&	9922	&	61/62		&	1823		&	5.4		&	28/40	&	1156	&	8.6	\\
Case 2	&	1.00	&	$\infty$	&	11/990	&	30532	&	54/55		&	1559		&	19.6		&	23/60	&	1662	&	18.4	\\
	&	0.10	&		&	18/1620	&	48666	&	173/174		&	4982		&	9.8		&	41/106	&	3100	&	15.7	\\
Case 3	&	1.00	&	$\infty$	&	16/12096	&	382828	&	202/203		&	7203		&	53.1		&	25/188	&	3367	&	113.7	\\
	&	0.10	&		&	3736/18144	&	563205	&	3736/3737	$^*$	&	422591	$^*$	&	1.3	\phantom{a.}$^*$	&	72/360	&	12375	&	45.5	\\
            \bottomrule
            \multicolumn{11}{r}{$\infty$: the model can not be solved by Gurobi. $*$: the algorithm cannot solve the problem to a 0.1\% tolerance but reach a 0.103\% tolerance}
        \end{tabular}
    }
\end{table}

From Table \ref{table:Comparative results for standard Benders, adaptive Benders and stabilised adaptive Benders.}, we can see that (a) stabilised Adaptive Benders is up to $113.7$ times faster than standard Benders for a 1.00\% convergence tolerance and $45.5$ times faster than standard Benders for a 0.10\% convergence tolerance, (b) Adaptive Benders gets slower when converging to a tighter tolerance and (c) compared with Adaptive Benders, stabilised Adaptive Benders is up to $2.14$ times faster for a 1.00\% convergence tolerance, and Adaptive Benders cannot solve the largest instance to 0.10\% due to severe oscillation. Therefore, for Case 3, we report the performance of unstabilised Adaptive Benders when it reaches a tolerance of 0.103\%, which is the tightest convergence tolerance it achieves and just before it starts oscillating severely.
\begin{table}[!htb]
    \centering
    \caption{Results of stabilised Adaptive Benders decomposition with different level sets.}
    \label{table:Results for stabilised adaptive Benders decomposition with different level sets.}
    \resizebox{\columnwidth}{!}{%
        \begin{tabular}
            {cS[table-format=1.2]|S[table-format=2.0]S[table-format=2.0]S[table-format=2.0]S[table-format=2.0]S[table-format=2.0]S[table-format=5.0]S[table-format=2.0]S[table-format=5.0]S[table-format=4.0]S[table-format=7.0]S[table-format=4.0]}
            \toprule
            \multicolumn{1}{c}{$\gamma$} & \multicolumn{1}{c}{$\epsilon$ (\%)} & \multicolumn{2}{c}{Case 0} & \multicolumn{2}{c}{Case 1} & \multicolumn{2}{c}{Case 2} & \multicolumn{2}{c}{Case 3} &\multicolumn{2}{c}{Average} \\ \cline{3-12}
            \multicolumn{1}{c}{} & \multicolumn{1}{c}{} & \multicolumn{1}{c}{Iters/Evals} & \multicolumn{1}{c}{Time (s)} & \multicolumn{1}{c}{Iters/Evals} & \multicolumn{1}{c}{Time (s)} & \multicolumn{1}{c}{Iters/Evals} & \multicolumn{1}{c}{Time (s)} & \multicolumn{1}{c}{Iters/Evals} & \multicolumn{1}{c}{Time (s)} & \multicolumn{1}{c}{Iters/Evals}& \multicolumn{1}{c}{Time (s)}\\ \hline
0.000*	&	1.00	&	26	/	35	&	970	&	21	/	45	&	1308	&	26	/	52	&	1396	&	25	/	188	&	5341	&	25	/	80	&	9015	\\
0.000	&	1.00	&	35	/	39	&	1106	&	21	/	45	&	1298	&	23	/	70	&	1925	&	24	/	169	&	7270	&	26	/	81	&	11599	\\
0.025	&	1.00	&	16	/	26	&	752	&	20	/	27	&	791	&	23	/	60	&	1662	&	30	/	105	&	3385	&	22	/	55	&	6590	\\
0.050	&	1.00	&	13	/	20	&	552	&	20	/	28	&	813	&	34	/	74	&	2039	&	21	/	163	&	4511	&	22	/	71	&	7915	\\
0.075	&	1.00	&	16	/	23	&	665	&	15	/	21	&	581	&	21	/	64	&	1859	&	30	/	232	&	10381	&	21	/	85	&	13486	\\
0.100	&	1.00	&	12	/	20	&	615	&	34	/	62	&	1820	&	17	/	60	&	1694	&	31	/	164	&	5419	&	24	/	77	&	9548	\\
0.200	&	1.00	&	18	/	24	&	657	&	35	/	64	&	1885	&	31	/	90	&	2588	&	46	/	180	&	12339	&	33	/	90	&	17469	\\
0.300	&	1.00	&	23	/	30	&	848	&	30	/	50	&	1478	&	33	/	88	&	2503	&	68	/	153	&	5055	&	39	/	80	&	9884	\\
0.400	&	1.00	&	22	/	29	&	842	&	37	/	58	&	1796	&	39	/	96	&	2755	&	140	/	633	&	22585	&	60	/	204	&	27978	\\
0.500	&	1.00	&	24	/	32	&	939	&	31	/	38	&	1276	&	58	/	143	&	4180	&	358	/	833	&	91031	&	118	/	262	&	97426	\\
0.600	&	1.00	&	24	/	31	&	980	&	44	/	57	&	1721	&	67	/	110	&	3457	&	566	/	1123	&	45590	&	175	/	330	&	51748	\\
0.700	&	1.00	&	32	/	37	&	1146	&	75	/	90	&	2970	&	126	/	212	&	6665	&	551	/	977	&	40753	&	196	/	329	&	51534	\\
0.800	&	1.00	&	39	/	45	&	1424	&	84	/	94	&	3438	&	197	/	275	&	9004	&	713	/	966	&	43529	&	258	/	345	&	57395	\\
0.900	&	1.00	&	73	/	81	&	2739	&	222	/	244	&	8165	&	487	/	649	&	21665	&	2149	/	2859	&	192864	&	733	/	958	&	225433	\\ \hline \hline
0.000*	&	0.10	&	51	/	61	&	1754	&	37	/	46	&	1343	&	45	/	128	&	3725	&	56	/	626	&	21070	&	47	/	215	&	27892	\\
0.000	&	0.10	&	55	/	59	&	1696	&	37	/	86	&	2508	&	41	/	116	&	3354	&	59	/	473	&	34296	&	48	/	184	&	41854	\\
0.025	&	0.10	&	35	/	47	&	1344	&	28	/	40	&	1156	&	41	/	106	&	3100	&	72	/	360	&	12472	&	44	/	138	&	18072	\\
0.050	&	0.10	&	21	/	33	&	941	&	29	/	43	&	1252	&	52	/	124	&	3693	&	62	/	502	&	16862	&	41	/	176	&	22748	\\
0.075	&	0.10	&	24	/	35	&	1022	&	27	/	37	&	1051	&	35	/	109	&	3181	&	81	/	690	&	49078	&	42	/	218	&	54332	\\
0.100	&	0.10	&	21	/	33	&	996	&	44	/	81	&	2384	&	33	/	105	&	3037	&	60	/	400	&	16605	&	40	/	155	&	23022	\\
0.200	&	0.10	&	26	/	36	&	1004	&	54	/	97	&	2880	&	54	/	153	&	4556	&	92	/	510	&	40627	&	57	/	199	&	49067	\\
0.300	&	0.10	&	31	/	43	&	1227	&	44	/	78	&	2322	&	61	/	88	&	7101	&	140	/	559	&	19737	&	69	/	192	&	30387	\\
0.400	&	0.10	&	27	/	37	&	1074	&	49	/	78	&	2395	&	53	/	174	&	5279	&	306	/	1499	&	59641	&	109	/	447	&	68389	\\
0.500	&	0.10	&	31	/	42	&	1234	&	39	/	50	&	1644	&	77	/	239	&	7181	&	471	/	1217	&	108217	&	155	/	387	&	118276	\\
0.600	&	0.10	&	31	/	41	&	1286	&	60	/	84	&	2513	&	112	/	211	&	6674	&	838	/	2054	&	93282	&	260	/	598	&	103755	\\
0.700	&	0.10	&	45	/	51	&	1555	&	89	/	115	&	3717	&	178	/	212	&	9965	&	954	/	2256	&	109567	&	317	/	659	&	124804	\\
0.800	&	0.10	&	50	/	62	&	1923	&	104	/	130	&	4523	&	326	/	599	&	18881	&	1085	/	1848	&	94686	&	391	/	660	&	120013	\\
0.900	&	0.10	&	85	/	94	&	3160	&	233	/	262	&	8709	&	657	/	1015	&	33098	&	2982	/	4706	&	376947	&	989	/	1519	&	421914	\\
            \bottomrule
            \multicolumn{12}{r}{$-$: algorithm stops because Gurobi fails solving a stabilisation problem. $*$: the stabilisation problem is removed}
        \end{tabular}
    }
\end{table}

\subsubsection{Improving the robustness}
The stabilisation factor $\gamma$ significantly impacts the performance. A very small $\gamma$ leads to loose stabilisation and makes stabilisation less effective, whereas a very large $\gamma$ leads to tight stabilisation and may hinder the exploitation of the solution space. We test the performance using different $\gamma$ from 0.025 to 0.9 and present the results in Table \ref{table:Results for stabilised adaptive Benders decomposition with different level sets.}. We find that a stabilisation factor less than 0.2 generally performs well. By checking the average performance over four cases, we find that  $\gamma$ equals 0.025 give the best performance. However, for different cases, the $\gamma$ that yields better performance varies. Furthermore, a rule of thumb for setting a fixed stabilisation factor may be to set it less than or equal to 0.5. 
% Although when $\gamma=0$, the target for \textbf{LMP} is the lower bound, it is still different from the unstabilised version. It is because the objective of \textbf{LMP} is to minimise the distance moved from the previous point while meeting the target. Therefore, when $\gamma=0$, \textbf{LMP} will get a point to reach the lower bound with minimum distance moved from the previous point.

\begin{table}[!htb]
    \centering
    \caption{Results for stabilised Benders decomposition with adjusted level sets (speed up: the time spent using fixed $\gamma$ divided by time spent using dynamic stabilisation).}
    \label{table:adjusting level sets case 0-3.}
        \resizebox{\columnwidth}{!}{%
    \begin{tabular}
        {cccS[table-format=1.2]S[table-format=0.1]|S[table-format=1.0]S[table-format=1.0]|S[table-format=5.0]S[table-format=1.2]S[table-format=1.0]S[table-format=1.0]S[table-format=1.2]}
        \toprule
        \multicolumn{1}{c}{} & \multicolumn{1}{c}{initial $\gamma$} & \multicolumn{1}{c}{$\omega$} & \multicolumn{1}{c}{$\underline{P}$} & \multicolumn{1}{c}{$\overline{P}$} & \multicolumn{2}{c}{Iters/Evals} & \multicolumn{4}{c}{Time (s)} \\ \cline{6-11} 
         \multicolumn{1}{c}{} & \multicolumn{1}{c}{} &\multicolumn{1}{c}{} & \multicolumn{1}{c}{}& \multicolumn{1}{c}{}& \multicolumn{1}{c}{$\epsilon=1.00\%$} & \multicolumn{1}{c}{$\epsilon=0.10\%$} &  \multicolumn{1}{c}{$\epsilon=1.00\%$} & \multicolumn{1}{c}{Speed up}&  \multicolumn{1}{c}{$\epsilon=0.10\%$} & \multicolumn{1}{c}{Speed up} \\ 
         \hline
Case 0	&	0.025	&	0.5	&	0.1	&	0.9	&	20	/	30	&	35	/	48	&	843	&	0.9	&	1368	&	1.0	\\
	&		&	0.9	&	0.1	&	0.9	&	17	/	27	&	28	/	44	&	702	&	1.1	&	1233	&	1.1	\\
	&	0.100	&	0.5	&	0.1	&	0.9	&	14	/	19	&	31	/	38	&	526	&	1.2	&	1081	&	0.9	\\
	&		&	0.9	&	0.1	&	0.9	&	13	/	18	&	23	/	31	&	484	&	1.3	&	862	&	1.2	\\
	&	0.500	&	0.5	&	0.1	&	0.9	&	17	/	24	&	27	/	37	&	753	&	1.3	&	1212	&	1.0	\\
	&		&	0.9	&	0.1	&	0.9	&	20	/	27	&	28	/	38	&	761	&	1.2	&	1112	&	1.1	\\
	&	0.900	&	0.5	&	0.1	&	0.9	&	31	/	39	&	53	/	43	&	1123	&	2.4	&	1554	&	2.0	\\
	&		&	0.9	&	0.1	&	0.9	&	33	/	41	&	43	/	55	&	1248	&	2.2	&	1655	&	1.9	\\ \hline
Case 1	&	0.025	&	0.5	&	0.1	&	0.9	&	22	/	31	&	35	/	49	&	882	&	0.9	&	1425	&	0.8	\\
	&		&	0.9	&	0.1	&	0.9	&	17	/	27	&	28	/	43	&	732	&	1.1	&	1238	&	0.9	\\
	&	0.100	&	0.5	&	0.1	&	0.9	&	14	/	27	&	31	/	57	&	753	&	2.4	&	1651	&	1.4	\\
	&		&	0.9	&	0.1	&	0.9	&	18	/	28	&	26	/	40	&	734	&	2.5	&	1157	&	2.1	\\
	&	0.500	&	0.5	&	0.1	&	0.9	&	22	/	42	&	31	/	54	&	1182	&	1.1	&	1599	&	1.0	\\
	&		&	0.9	&	0.1	&	0.9	&	22	/	45	&	34	/	45	&	1254	&	1.0	&	2054	&	0.8	\\
	&	0.900	&	0.5	&	0.1	&	0.9	&	15	/	20	&	30	/	44	&	627	&	13.0	&	1342	&	6.5	\\
	&		&	0.9	&	0.1	&	0.9	&	26	/	41	&	35	/	61	&	1277	&	6.4	&	1901	&	4.6	\\ \hline
Case 2	&	0.025	&	0.5	&	0.1	&	0.9	&	38	/	68	&	51	/	110	&	1979	&	0.8	&	3224	&	1.0	\\
	&		&	0.9	&	0.1	&	0.9	&	18	/	61	&	41	/	152	&	1601	&	1.0	&	4461	&	0.7	\\
	&	0.100	&	0.5	&	0.1	&	0.9	&	28	/	67	&	36	/	94	&	1775	&	1.0	&	2720	&	1.1	\\
	&		&	0.9	&	0.1	&	0.9	&	29	/	67	&	39	/	99	&	1931	&	0.9	&	2858	&	1.1	\\
	&	0.500	&	0.5	&	0.1	&	0.9	&	39	/	90	&	46	/	108	&	2533	&	1.7	&	3191	&	2.3	\\
	&		&	0.9	&	0.1	&	0.9	&	39	/	102	&	49	/	141	&	3124	&	1.3	&	4261	&	1.7	\\
	&	0.900	&	0.5	&	0.1	&	0.9	&	48	/	95	&	76	/	194	&	2703	&	8.0	&	5839	&	5.7	\\
	&		&	0.9	&	0.1	&	0.9	&	50	/	107	&	65	/	164	&	3377	&	6.4	&	5086	&	6.5	\\ \hline
Case 3	&	0.025	&	0.5	&	0.1	&	0.9	&	49	/	217	&	84	/	418	&	6573	&	0.7	&	14067	&	1.2	\\ 
	&		&	0.9	&	0.1	&	0.9	&	28	/	173	&	66	/	399	&	5486	&	0.8	&	15195	&	1.1	\\ 
	&	0.100	&	0.5	&	0.1	&	0.9	&	24	/	111	&	87	/	433	&	3306	&	2.0	&	14860	&	1.0	\\ 
	&		&	0.9	&	0.1	&	0.9	&	44	/	143	&	126	/	606	&	4538	&	1.5	&	21276	&	0.7	\\
	&	0.500	&	0.5	&	0.1	&	0.9	&	52	/	182	&	140	/	566	&	5837	&	7.5	&	24892	&	2.9	\\
	&		&	0.9	&	0.1	&	0.9	&	61	/	210	&	123	/	494	&	6901	&	6.4	&	17332	&	4.2	\\
	&	0.900	&	0.5	&	0.1	&	0.9	&	58	/	194	&	98	/	425	&	5758	&	33.5	&	14838	&	25.4	\\
	&		&	0.9	&	0.1	&	0.9	&	76	/	178	&	133	/	444	&	6042	&	31.9	&	15820	&	23.8	\\ \hline
\multicolumn{4}{c}{Average}							&		&	31	/	80	&	56	/	174	&	2417	&		&	6011	&		\\
\multicolumn{4}{c}{Average (fixed $\gamma$)$^*$}							&		&	225	/	338	&	307	/	550	&	84749	&		&	145321	&		\\
\multicolumn{4}{c}{Standard deviation}							&		&	16	/	64	&	35	/	183	&	2052	&		&	6854	&		\\
\multicolumn{4}{c}{Standard deviation (fixed $\gamma$)$^*$}							&		&	342	/	424	&	458	/	656	&	102821	&		&	190074	&		\\
        \bottomrule
        \multicolumn{11}{r}{$*$: consider the runs when $\gamma=0.025,\ 0.1,\ 0.5 \text{ and } 0.9$ }
    \end{tabular}
    }
\end{table}

We see that different fixed $\gamma$ can lead to a noticeable difference in performance. It is aligned with Remark 5 in \citep{Zverovich2012}. Unlike \citep{Zverovich2012} who decided to use $\gamma=0.5$ to set the level set, we test the approach presented in Section \ref{sec:algorithm} to adjust $\gamma$ and then the level set that may make the stabilisation more robust and independent of the choice of the stabilisation factor. 

We test extensively the dynamic stabilisation scheme on Case 0-3. The results for different cases are shown in Table \ref{table:adjusting level sets case 0-3.}. By comparing results from Table \ref{table:adjusting level sets case 0-3.} and their fixed $\gamma$ counterparts in Table \ref{table:Results for stabilised adaptive Benders decomposition with different level sets.}, we can see that in almost all cases adjusting the level set outperforms the cases using a fixed $\gamma$. Although sometimes dynamic stabilisation is slightly slower than standard stabilisation, most of the time it is better. Furthermore, dynamic stabilisation makes the level set method stabilisation more robust in terms of the choice of initial $\gamma$. It is particularly valuable because one may need extensive tests to find the $\gamma$ that yields the best performance for the problems to be solved. However, as we see in Table \ref{table:Results for stabilised adaptive Benders decomposition with different level sets.}, different problems may have different best $\gamma$. Therefore, a dynamic stabilisation that makes the performance less dependent on the choice of $\gamma$ may make it easier to get a satisfying performance if one chooses a bad initial $\gamma$ because the dynamic adjustment will help correct $\gamma$ to a sensible value while solving the problem.

\subsection{Power system analysis}
\label{sec:power_system_analysis}
In this section, we present the results of the 5-region UK power system planning problem. We analyse the investment decisions, expected costs, and the Value of the Stochastic Solution (VSS). 

The investment decisions in the first stage are presented in Table \ref{table:Investment in onshore wind in the first investment stage.}. There are no investments in technologies except onshore wind in the first investment stage. We notice that the transmission lines are expanded in the later investment nodes. Therefore, for the first investment stage, only investment in onshore wind is presented. The onshore wind is mainly invested in Scotland, North England and South England in the first investment stage. When considering only short-term uncertainty, we can see that in Case 0, a total of 90.85 GW of onshore wind is invested, 28\% of which is in North England. Compared with Case 1, around 3.8 GW less capacity is installed in Case 0. When considering uncertainty in both long-term demand and CO$_2$ budget, we can see a 3.56 GW investment in onshore wind in Scotland, compared with 14.76 GW in Case 1 and 1.42 GW in Case 3. 

\begin{table}[!htb]
    \centering
    \caption{\centering Investment in onshore wind in the first investment stage.}
    \label{table:Investment in onshore wind in the first investment stage.}
    % \vspace{-0.45cm}
    \resizebox{\columnwidth}{!}{%
        \begin{tabular}
            {c|rrrrrr}
            \toprule
            \multicolumn{1}{c}{} & \multicolumn{6}{c}{Investment in onshore wind (GW)}  \\ \cline{2-7}
            \multicolumn{1}{c}{}   & \multicolumn{1}{c}{Scotland} & \multicolumn{1}{c}{North England} & \multicolumn{1}{c}{Midlands \& Wales} & \multicolumn{1}{c}{East England} & \multicolumn{1}{c}{South England}  & \multicolumn{1}{c}{UK Total}\\ \hline
            Case 0	&	0.00	&	25.70	&	0.00	&	0.00	&	65.15	&	90.85	\\
            Case 1	&	14.76	&	28.25	&	0.00	&	0.00	&	51.68	&	94.69	\\
            Case 2	&	3.56	&	23.22	&	0.00	&	0.00	&	59.19	&	85.97	\\
            Case 3	&	1.42	&	23.60	&	0.00	&	0.00	&	69.16	&	94.18	\\
            \bottomrule
    \end{tabular}
    }
\end{table}

Table \ref{table:expected costs and VSS.} shows the optimal costs and the VSS for considering long-term uncertainties. We can see that there is up to £$7,702$ million VSS when considering uncertainty, including CO$_2$ budget and long-term demand. The VSS is £$2,904$ million when considering only CO$_2$ budget as an uncertainty parameter. When considering long-term uncertainty, including CO$_2$ budget, CO$_2$ tax and long-term demand, the VSS is 4.4\% of the optimal cost. This shows the value of including long-term uncertainty in a long-term planning problem and solving a large model. 
\begin{table}[!htb]
    \centering
    \caption{\centering Optimal costs and VSS.}
    \label{table:expected costs and VSS.}
    % \vspace{-0.45cm}
    % \resizebox{\columnwidth}{!}{%
        \begin{tabular}
            {c|rrrr}
            \toprule
            \multicolumn{1}{c}{}   & \multicolumn{1}{c}{Case 0} & \multicolumn{1}{c}{Case 1} & \multicolumn{1}{c}{Case 2} & \multicolumn{1}{c}{Case 3} \\ \hline
            Optimal cost (mn £) &174099  &174276  &174871 &174785\\
            VSS (long-term uncertainty, mn £) &-  &2904 &7325  &7702 \\
            \bottomrule
    \end{tabular}
    % }
\end{table}

\section{Discussion}
\label{sec:discussion}
In this paper, we propose a method to address the computational difficulty of a multi-stage stochastic programming problem with short-term and long-term uncertainty that is formulated using a multi-horizon stochastic programming approach. Similar studies on developing a method to solve such type of problems can be found in \citep{Zakeri, Downward2020}. Compared with their approach, we exploit the properties of the \textbf{SP} and stabilise the algorithm with the level set method and adaptively select \textbf{SP}s to solve exactly per iteration for better approximation, which shows significant performance improvement. The method can be generally applied to solve any problem that is formulated in \eqref{eq:MP} and \eqref{eq:SP}. 

We demonstrate our proposed method on a multi-region UK power system planning problem. To the authors' knowledge, this is the first study that presents a multi-horizon formulation of a multi-region power system planning problem with short-term and long-term uncertainty and proposes a method to solve such a problem efficiently. Compared with a similar problem for long-term investment planning such as \citep{Backe2022} that only considers short-term uncertainty, this paper firstly introduces both long-term and short-term uncertainty in a power system planning problem using a multi-horizon framework. 

We notice and analyse the oscillation of the Benders-type decomposition method for multi-region investment planning problems. The level set method stabilisation approach was used to stabilise Benders. Compared with the existing literature that studied the level set method, we integrate it with the inexact oracles and show that it significantly improves computational performance. In addition, similar studies normally set the level set in an ad hoc way \citep{Zverovich2012,Ruszczynski1997}. Moreover, we test to adjust the target based on a proposed measurement. For the test instance, adjusting the level set can usually yield better or equivalent performance. However, the parameters that yield the best performance may be case-dependent.

Although the stabilisation is useful, the stabilisation problem can potentially be a large QP and slow to solve. One possible approach to stabilise the problem efficiently is to utilise the built-in method, analytic centre \citep{Gondzio1996} in a commercial solver like Gurobi to solve a feasibility problem to avoid solving a QP. We test utilising the analytic centre of Gurobi to potentially avoid solving a QP \textbf{LMP}. However, the results show that proper stabilisation may still be the better option, even for large problems.

We demonstrate the method for solving large-scale linear programming. However, the method can be applied to solve mixed-integer linear programming problems without modification as long as the integer variables are in the \textbf{MP}. In such a case, the stabilisation problem becomes a mixed-integer QP problem which may be slow to solve. Some other stabilisation techniques, such as local branching \citep{Baena2020a} may be an alternative. 

Although this paper presents a general method to solve a class of large-scale optimisation problems very efficiently, there are some limitations. Firstly, we need the same coefficient matrices in all nodes to utilise the inexact oracles. This may be limited when different operational scenarios are preferred. However, the inexact oracles can be generalised to apply to some groups of nodes with the same matrices and work in a problem with different matrices in different nodes. However, having different scenarios in different nodes may lead to lower stability of a stochastic problem. Secondly, although multi-horizon formulation significantly reduces the problem size, it may be limited when long-term storage, such as pumped hydro storage, is in the system. However, an easy fix to this issue is to collect information on the storage level at the end of one stage in the master problem and pass it to the next stage. Thirdly, in the case study, we only demonstrate the proposed method to solve the UK power system planning problem and show significant performance improvement, but more problems may be solved using the proposed method to gain more insights into the performance of the algorithm.

\section{Conclusions and future work}
\label{sec:conclusions}
In this paper, we proposed stabilised Benders decomposition with adaptive oracles to solve long-term multi-region investment planning problems with short-term and long-term uncertainty. We applied the algorithm to solve a multi-region UK power system investment planning problem towards 2035. We formulated such a problem using a multi-horizon stochastic programming approach. The test instances have up to $1$ billion variables and $4.5$ billion constraints. The computational results show that: a) for a $1.00\%$ convergence tolerance, the proposed stabilised method is up to $113.7$ times faster than standard Benders decomposition and $2.14$ times faster than Adaptive Benders decomposition without stabilisation; b) for a $0.10\%$ convergence tolerance, the proposed stabilised method is up to $45.5$ times faster than standard Benders decomposition and the unstabilised Adaptive Benders decomposition cannot solve the largest instance to the convergence tolerance due to severe oscillation and c) dynamic level set method increases the robustness of the stabilisation.

Although the proposed method reduced the computational effort significantly and was used to solve multi-horizon stochastic programming with short-term and long-term uncertainty, we notice that for a very large problem with many decision nodes, the reduced master problem and the stabilisation problem may take longer to solve. Therefore, in future, techniques including node aggregation and cuts selection and deletion may be needed to improve the performance. In addition, although multi-horizon formulation reduces the problem size significantly, the model size may be reduced further by adjusting the scenario tree, e.g., removing the scenarios that do not make a difference while solving the problem. 

\section*{CRediT author statement}
\textbf{Hongyu Zhang:} Conceptualisation, Methodology, Software, Validation, Formal analysis, Investigation, Visualisation, Data curation, Writing - original draft, Writing - review \& editing. \textbf{Nicolò Mazzi:} Conceptualisation, Methodology, Data curation, Software, Writing - review \& editing. \textbf{Ken McKinnon:} Conceptualisation, Methodology, Supervision, Writing - review \& editing, Funding acquisition. \textbf{Rodrigo Garcia Nava:} Conceptualisation, Methodology, Software, Writing - review \& editing. \textbf{Asgeir Tomasgard:} Conceptualisation, Supervision, Writing - review \& editing, Funding acquisition.

\section*{Declaration of competing interest}
The authors declare that they have no known competing financial interests or personal relationships that could have appeared to influence the work reported in this paper.

\section*{Acknowledgements}
This work was supported by the Research Council of Norway through PETROSENTER LowEmission [project code 296207]; and the Engineering and Physical Sciences Research Council (EPSRC) through the CESI project [EP/P001173/1].

\appendix
\section{Appendix}
\label{nomenclature}
% \printnomenclature[0.9in]
\begin{description}[itemsep=-5pt, leftmargin=!,labelwidth=\widthof{\bfseries $p_{st}^{SE+}/p_{st}^{SE-}$}]
    \item [\textbf{Investment planning model sets}]
    \item [{$\mathcal{P}$}]set of technologies
    \item [{$\mathcal{I}$}]set of operational nodes, indexed by $i$
    \item [{$\mathcal{I}_{0}$}]set of investment nodes, indexed by $i_{0}$
    \item [{$\mathcal{I}_{i}$}]set of investment nodes $i$ $(i \in \mathcal{I}_{0})$ ancestor to operational node $i$ $(i \in \mathcal{I})$
    \item [\textbf{Operational model sets}]
    \item [{$\mathcal{N}$}]set of time slices
    \item [{$\mathcal{T}$}]set of hours in all time slices
    \item [{$\mathcal{L}$}]set of transmission lines
    \item [{$\mathcal{G}$}]set of thermal generators
    \item [{$\mathcal{S}$}]set of electricity storage
    \item [{$\mathcal{R}$}]set of renewable generations
    \item [\textbf{Investment planning model parameters}]
    \item [{$C^{Inv}_{pi}$}]unitary investment cost of device $p$ in investment node $i$ ($p \in \mathcal{P}, i \in \mathcal{I}_{0}$) [£/MW]
    \item [{$C^{Fix}_{p}$}]unitary fix operational and maintenance cost of device $p$ ($p \in \mathcal{P}$) [£/MW]
    \item [{$X^{Hist}_{p}$}]historical capacity of device $p$ ($p \in \mathcal{P}$) [MW]
    \item [{$X^{Max}_{p}$}]maximum installed capacity of device $p$ ($p \in \mathcal{P}$) [MW]
    \item [{$\kappa$}]scaling effect depending on the number of operation years between investment nodes
    \item [{$\delta^{I_{0}}_{i}/\delta^{I}_{i}$}]discount factor of investment node $i$ ($i_{0} \in \mathcal{I}_{0}$)/ operational node $i$ ($i \in \mathcal{I}$)
    \item [{$\pi^{I_{0}}_{i}/\pi^{I}_{i}$}]probability of investment node $i$ ($i_{0} \in \mathcal{I}_{0}$)/ operational node $i$ ($i \in \mathcal{I}$)
    \item [{$H^{P}_{p}$}]life time of technology $p$ ($p \in \mathcal{P}$)
    \item [{$x_{i}$}]right hand side coefficients of the operational subproblem
    \item [{$c_{i}$}]cost coefficients of the operational subproblem
    \item [{$\mu^{E}_{i}$}]CO$_2$ budget at operational node $i$ ($i \in \mathcal{I}$)
    \item [{$\mu^{DP}_{i}$}]scaling factor on power demand at operational node $i$ ($i \in \mathcal{I}$)
    \item [{$C^{CO2}_{i}$}]CO$_2$ emission price at operational node $i$ ($i \in \mathcal{I}$)
    \item [\textbf{Operational model parameters}]
    \item [{$\mu^{E}$}]yearly CO$_2$ emission limit (tonne)
    \item [{$\pi_{t}$}]probability of operation period $t$ ($t \in \mathcal{T}$)
    \item [{$H_t$}]number of hour(s) in one operational period $t$
    \item [{$\alpha^{G}_{g}$}]maximum ramp rate of gas turbines ($g \in \mathcal{G}$) [MW/MW]
    \item [{$R_{rt}^{R}$}]capacity factor of renewable unit $r$ in period $t$ ($r \in \mathcal{R}, t \in \mathcal{T}$)
    \item [{$\eta_{s}^{SE}$}]efficiency of electricity store $s$ ($s \in \mathcal{S}$)
    \item [{$\gamma_{s}^{SE}$}]power ratio of electricity store $s$ ($s \in \mathcal{S}$) [MWh/MW]
    \item [{$E^{G}_g$}]emission factor of gas turbine $g$ $(g \in \mathcal{G})$ [tonne/MWh]
    \item [{$C_g^G/C_s^{SE}$}]total operational cost of a generator $g$/ a storage facility $s$ ($g \in \mathcal{G}$/ $s \in \mathcal{S}$) [£/MW]
    \item [{$C^{ShedP}$}]power load shed penalty cost [£/MW]
    \item [{$P^{DP}_{zt}$}]power demand at region $z$ period $t$ $(z \in \mathcal{Z}, t \in \mathcal{T})$ [MW]
    \item [\textbf{Investment planning model variables}]
    \item [{$x_{pi}^{Acc}$}]accumulated capacity of device $p$ in operational node $i$ ($p \in \mathcal{P}, i \in \mathcal{I}$) [MW]
    \item [{$x_{pi}^{Inst}$}]newly invested capacity of device $p$ in investment node $i_{0}$ ($p \in \mathcal{P}, i \in \mathcal{I}_{0}$) [MW]
    \item [{$c_{i}^{INV}$}]total expected investment cost (£)
    \item [\textbf{Operational model variables}]
    \item [{$p_{g}^{AccG}$}]accumulated capacity of gas turbine $g$ $(g \in \mathcal{G})$ [MW]
    \item [{$p_{g}^{AccG}$}]accumulated capacity of renewable unit $r$ $(r \in \mathcal{R})$ [MW]
    \item [{$p_{s}^{AccSE}$}]accumulated charging/discharging capacity of electricity store $s$ ($s \in \mathcal{S}$) [MW]
    \item [{$p_{gt}^{G}$}]power generation of gas turbine $g$ in period $t$ ($g \in \mathcal{G}, t \in \mathcal{T}$) [MW]
    \item [{$p_{st}^{SE+}/p_{st}^{SE-}$}]charge/discharge power of electricity store $s$ in period $t$ ($s \in \mathcal{S}, t \in \mathcal{T}$) [MW]
    \item [{$p_{lt}^{L}$}]power flow in line $l$ in period $t$ ($l \in \mathcal{L}, t \in \mathcal{T}$) [MW]
    \item [{$p_{zt}^{GShedP}$}]generation shed at $z$ in period $t$ ($z \in \mathcal{Z}, t \in \mathcal{T}$) [MW]
    \item [{$q_{st}^{SE}$}]energy level of electricity store $s$ at the start of period $t$ ($s \in \mathcal{S}, t \in \mathcal{T}$) [MWh]
    \item [{$p_{zt}^{ShedP}$}]load shed at $z$ in period $t$ ($z \in \mathcal{Z}, t \in \mathcal{T}$) [MW]
    \item [\textbf{Function}]    
    \item [{$c^{OPE}(\cdot,\cdot)$}]operational cost at node at operational node $i$ ($i \in \mathcal{I}$) (£)

\end{description}



 % \footnote{For further information, please refer to \url{https://docs.julialang.org/en/v1/manual/asynchronous-programming/}}

% \section*{References}
%\bibliographystyle{elsarticle-harv}
\clearpage
\setlength{\bibsep}{0pt plus 0.3ex}
\footnotesize{
\bibliographystyle{model5-names}
\bibliography{benders_ao_lm_v4}}
\end{document}